\documentclass[preprint,authoryear,12pt]{elsarticle}
\usepackage[sumlimits]{amsmath}
\usepackage{amsfonts}
\usepackage{color}
\usepackage{graphicx}
\usepackage{xypic}
\journal{Ocean Modeling}

\DeclareMathOperator{\on}{on}

\DeclareMathOperator{\diff}{d}

\DeclareMathOperator{\Top}{\text{top}}
\DeclareMathOperator{\Floor}{\text{floor}}
\DeclareMathOperator{\Coast}{\text{coast}}

\DeclareMathOperator{\Cond}{Cond}

\def\MM#1{\boldsymbol{#1}}
\newcommand{\pp}[2]{\frac{\partial #1}{\partial #2}}

\newcommand{\eqnref}[1]{(\ref{#1})}

\begin{document}
\begin{frontmatter}
\title{Solving the Poisson equation on small aspect ratio domains
  using unstructured meshes} 

\author[ese]{S. C. Kramer\corref{cor1}}
\cortext[cor1]{s.c.kramer@imperial.ac.uk}
\author[aero]{C. J. Cotter}
\author[ese]{C. C. Pain}
\address[ese]{Department of Earth Science and Engineering,
    Imperial College London, Prince Consort Road, London SW7 2AZ, UK}
\address[aero]{Department of Aeronautics, Imperial College
    London, Prince Consort Road, London SW7 2AZ, UK}

\begin{abstract} We discuss the ill conditioning of the matrix for the
  discretised Poisson equation in the small aspect ratio limit, and
  motivate this problem in the context of nonhydrostatic ocean
  modelling. Efficient iterative solvers for the Poisson equation in
  small aspect ratio domains are crucial for the successful
  development of nonhydrostatic ocean models on unstructured meshes.
  We introduce a new multigrid preconditioner for the Poisson problem
  which can be used with finite element discretisations on general
  unstructured meshes; this preconditioner is motivated by the fact
  that the Poisson problem has a condition number which is independent
  of aspect ratio when Dirichlet boundary conditions are imposed on
  the top surface of the domain. This leads to the first level in an
  algebraic multigrid solver (which can be extended by further
  conventional algebraic multigrid stages), and an additive smoother. We
  illustrate the method with numerical tests on unstructured meshes,
  which show that the preconditioner makes a dramatic improvement on a
  more standard multigrid preconditioner approach, and also show that
  the additive smoother produces better results than standard SOR
  smoothing. This new solver method makes it feasible to run
  nonhydrostatic unstructured mesh ocean models in small aspect ratio
  domains.
\end{abstract}


\end{frontmatter}
\section{Introduction}

There are many processes in the ocean (such as separating Western
boundary currents, and density overflows) which are small in scale and
restricted to particular regions, but which form crucial components of
the global ocean circulation mechanism. It therefore seems attractive
to design ocean models which use the finite element method on fully
unstructured meshes in order to incorporate some of these smaller
scale features into a global ocean model (see \citet{Pain_etal2005}
for background and references). However, there are numerous pitfalls
to negotiate in order to achieve this goal, arising from the fact that
the global ocean is very thin: the horizontal lengthscale is thousands
of times larger than the vertical lengthscale.

One particular issue arises
if one wishes to relax the hydrostatic approximation
\citep{Pe1987}, 
allowing for a model which is valid on both small and
large scales. The nonhydrostatic pressure is obtained by solving a
three dimensional elliptic problem with very large eigenvalues
resulting from the horizontal scales, as well as very small
eigenvalues resulting from the vertical scales. This means that the
system is very ill conditioned. Since the Conjugate Gradient method,
which is typically used for finite element discretisations with many
degrees of freedom, has a convergence rate which scales with the
square root of the condition number (see \citet{Sh1994} for example),
this can have a catastrophic effect on the performance of the
numerical model.

In the ocean modelling context this problem was first
encountered by \citet{marshall97}. A solution strategy was proposed
using a vertical preconditioner which solves the vertically integrated
(aspect ratio independent) equations and then distributes the solution
throughout the mesh. It was shown that the use of this preconditioner
resulted in nonhydrostatic simulations which were as fast as
hydrostatic simulations at the same resolution. This strategy, has
since been used in a number of nonhydrostatic ocean models, including
those on horizontally unstructured grids such as
\citet{FrGeSt2006}. However, the vertical averaging depends on the
computational mesh being organised in vertical layers. This prohibits
more general types of vertically unstructured meshes which may be
required for multiscale simulations in which a small scale process is
resolved within a large scale flow, or for hybrid meshes which
accommodate both terrain following and isopycnal (constant density)
layers.

In this paper we extend the vertical averaging strategy so that it can
be applied to vertically unstructured meshes of large-scale ocean
modelling such as those that can be used in the Imperial College Ocean
Model (ICOM) \citep{Pi_etal2008}. The extension is formulated by using
the vertical extrapolation operator, which takes any point in the
domain and returns the value of a function at the top surface directly
above that point. This operator is the dual of the vertical
integration operator, and can easily be approximated on a vertically
unstructured mesh. This extension is described within the context of
the algebraic multigrid method; the framework also shows how to
incorporate further algebraic multigrid stages into the ``smoother'',
which reconstructs the solution of the vertically integrated equations
throughout the domain. This turns out to be necessary when a genuinely
multiscale mesh is used in which the aspect ratio becomes
$\mathcal{O}(1)$ at the smallest scales. The aim is to obtain a
numerical solver which has a convergence rate which is independent of
the aspect ratio.

The rest of this paper is organised as follows. In section
\ref{background}, we describe the type of problems we wish to solve on
small aspect ratio domains, and motivate them using the ocean
modelling applications. In particular, this section explains why we
cannot avoid solving an elliptic problem with Neumann boundary
conditions on all surfaces, which is precisely the problem which gives
rise to ill-conditioning. In section \ref{finite element} we formulate
the problem as a finite element approximation, in order to fix
notation, and in section \ref{condition number section} we compute some
estimates on the condition number for the Neumann boundary condition
case, as well as the case where Dirichlet boundary conditions are
imposed on the top surface. It is observed that imposing the Dirichlet
boundary conditions removes the small eigenvalues, and this motivates
a preconditioning strategy in which one first eliminates the interior
degrees of freedom to obtain an equation for the solution on the top
surface, then one uses this surface solution as a Dirichlet boundary
condition to reconstruct the solution throughout the domain. This
paves the way for section \ref{new preconditioner} in which our
proposed preconditioner is introduced, in the context of algebraic
multigrid preconditioners for the Conjugate Gradient method. The
preconditioner is tested in various examples in section \ref{numerics}.
Finally, in section \ref{summary} we give a summary and outlook.
 
\section{Background: oceanographic applications}
\label{background}
In this section we describe how the pressure Poisson equation arises
in nonhydrostatic models, which motivates the need to develop
efficient solvers for this equation in small aspect ratio domains. We
shall also explain the types of boundary conditions that are imposed,
in particular we shall explain why we need to tackle the problem of
solving the pressure Poisson equation with Neumann boundary
conditions. 

\subsection{Nonhydrostatic equations}
The nonhydrostatic Euler-Boussinesq equations for a rotating
stratified fluid on an $f$-plane are
\begin{eqnarray}
\label{momentum}
\rho_0\left(\pp{\MM{u}}{t} + \MM{u}\cdot\nabla\MM{u}
+ f\hat{\MM{z}}\times\MM{u}\right) & = & -\nabla p - g\rho\hat{\MM{z}}, \\
\label{incompressibility}
\nabla\cdot\MM{u} & = & 0, \\
\pp{T}{t} + \MM{u}\cdot\nabla T & = & 0, 
\end{eqnarray}
where $\MM{u}$ is the velocity, $\rho_0$ is the (constant) reference
density, $\hat{\MM{z}}$ is the unit vector in the upward direction, $f$ is
the Coriolis parameter, $p$ is the pressure, $g$ is the gravitation
constant, $T$ is the potential temperature and $\rho$ is a prescribed
function of the temperature.  Here we do not complicate the exposition
by including viscosity or diffusivity terms, or dependence of the
density on salinity or pressure. 

We solve the problem in an ``ocean shaped'' domain $\Omega$ with
bottom boundary $\partial\Omega_{\Floor}$, coastal boundaries
$\partial\Omega_{\Coast}$ and top boundary $\partial\Omega_{\Top}$ (which may
be allowed to move up and down to accommodate surface waves). The
horizontal extent of the domain is $L$, and the vertical extent of the
domain is $H$; in this paper we concentrate on the difficulties when
the aspect ratio $H/L$ of the domain is very small ($H/L\approx 1/1000$
for an ocean basin). We shall parameterise the top surface by
\begin{equation*}
z = \eta(x,y),
\end{equation*}
with $\eta=0$ when the fluid is at rest, and make the additional
simplifying assumption that the coastal boundaries are vertical.

 We consider two types of boundary conditions:
\begin{itemize}
\item {\bfseries Rigid lid}: 
\begin{eqnarray*}
\MM{u}\cdot\MM{n}&=&0, \quad \MM{x} \in \partial\Omega_{\Floor}\cup
\partial\Omega_{\Coast}\cup \partial\Omega_{\Top}, \\
 \eta &= & 0,
\end{eqnarray*}
for a constant top surface height $z=0$.
\item {\bfseries Free surface}: 
\begin{eqnarray}
\nonumber \MM{u}\cdot\MM{n}&=&0, \quad \MM{x} \in \partial\Omega_{\Floor}\cup
\partial\Omega_{\Coast}, \\
\nonumber p & = & p_a, \quad \MM{x} \in \partial\Omega_{\Top}, \\
 \eta_t & = & - \MM{u}_H\cdot\nabla_H\eta-w
= -\frac{\MM{u}\cdot\MM{n}}{\MM{n}\cdot\hat{\MM{z}}}.
\label{kinematic}
\end{eqnarray}
\end{itemize}
In the subsequent section we shall show that both of these equations
effectively result in an ill conditioned pressure equation in the
small aspect ratio limit; this always occurs for the rigid lid case
and also occurs for the free surface case when one wishes to take
large timesteps and has an unstructured mesh.

\subsection{Solving for the pressure}
The role of the pressure in these equations is as a Lagrange
multiplier which enforces the incompressibility condition
\eqref{incompressibility} (in fact it enforces that the solution
to $D_t+ \nabla\cdot(D\MM{u})=0$ is $D=1$, and equation
\eqref{incompressibility} is then a direct consequence, see 
\citet{HoMaRa1998}) for example. To solve for the pressure, take the divergence
of equation \eqref{momentum}:
\begin{equation*}
-\nabla^2 p =
\nabla\cdot\left(
\rho_0\left(\MM{u}\cdot\nabla + f\hat{\MM{z}}\times\right)\MM{u} + 
\hat{\MM{z}}g\rho\right),
\end{equation*}
which is a Poisson equation for the pressure $p$, given the other
variables $\MM{u}$ and $T$. In practise, one usually solves for the
pressure update $\rho_0\phi/\Delta t = p^{n+1}-p^n$ (where $p^n$ is
the pressure at time level $n$) using one of the families of
projection methods based on \cite{Ch1967,Te1969}. These methods are
typically predictor-corrector schemes in which a predictor $\MM{u}^*$
is obtained using the pressure from the previous timestep, without
enforcing the incompressibility condition \eqref{incompressibility},
and then constructing a correction
\begin{equation}
  \label{update}
  \MM{u}^{n+1} = \MM{u}^* - \nabla\phi,
\end{equation}
subject to
\begin{equation*}
\nabla\cdot\MM{u}^{n+1} = 0.
\end{equation*}
Taking the divergence of equation \eqref{update} gives
\begin{equation}
  \label{phi eqn}
  -\nabla^2\phi = -\nabla\cdot\MM{u}^*.
\end{equation}
For more details, see \cite{KaSh2005,GrSa2000}.

On slip boundaries where $\MM{u}\cdot\MM{n}=0$, since $\MM{u}^*$
already satisfies the boundary conditions, then 
\begin{equation}
  \label{neumann phi}
  \pp{\phi}{n} = 0,
\end{equation}
so that $\MM{u}^{n+1}$ keeps the same boundary conditions. In the free
surface case, the pressure must stay constant on the surface, so
$\phi$ satisfies the zero Dirichlet condition on the top surface:
\begin{equation*}
  \pp{\phi}{n}=0 \text{ on }
  \partial\Omega_{\Coast}\cup\partial\Omega_{\Floor}, \quad
  \phi = 0 \text{ on }\partial\Omega_{\Top}.
\end{equation*}
However, the equations support free surface (barotropic) waves which
are very fast compared to the other dynamics. Moreover, the highest
frequency waves have small amplitudes and we are not concerned with
the details of their evolution in large scale models. The standard
ocean modelling approach is to apply a splitting method in which one
integrates the barotropic waves in time using a much smaller timestep;
the remaining (baroclinic) dynamics is then integrated using a larger
timestep (see \citet{ShMc2005} for example).  However, if the mesh
used is unstructured in the vertical direction (\emph{i.e.} the mesh
is not arranged into layers) then it is not possible to obtain such a
splitting. An alternative approach which is often used in coastal
engineering applications (see \citet{LaPi2005} for example), is to
construct a new ``piezometric'' pressure
\begin{equation*}
  \tilde{p}(x,y,z) = p(x,y,z) - \rho_0g\eta(x,y) - p_a,
\end{equation*}
which specifies the Dirichlet boundary condition
\begin{equation}
  \label{piezo bc}
  \tilde{p} = -\rho_0g\eta, \quad \on\,\,\partial\Omega_{\Top}.
\end{equation}
This means that we do not need a separate free surface variable as the
value can simply be read from $\tilde{p}$ at the free surface.  This
piezometric variable satisfies the same pressure Poisson equation with
modified right-hand side and different boundary conditions.

Taking the time derivative of equation \eqref{piezo bc} and
substituting the kinematic boundary condition \eqref{kinematic} gives
\begin{equation*}
  \pp{\tilde{p}}{t} = \rho_0g\frac{\MM{u}\cdot\MM{n}}{\MM{n}\cdot\MM{k}},
    \quad \on\,\,\partial\Omega_{\Top}.
\end{equation*}
If we wish to take large timesteps compared to the timescale of the
free surface waves, the barotropic-baroclinic split is not available
on unstructured meshes and it is necessary to solve the equations
using a linearly implicit method. Here we will describe the simplest
case, the backward Euler scheme, but the development of higher order
schemes is similar. The backward Euler scheme gives
\begin{equation*}
  \frac{p^{n+1}-p^n}{\Delta t} =
    \rho_0g\frac{\MM{u^{n+1}}\cdot\MM{n}}{\MM{n}\cdot\MM{k}}, \quad
    \on\,\,\partial\Omega_{\Top}.
\end{equation*}
Substituting the pressure update gives
\begin{equation*}
  \frac{\phi}{\Delta t^2} = 
  g\frac{\pp{\phi}{n}}{\MM{n}\cdot\MM{k}}
  , \quad
  \on\,\,\partial\Omega_{\Top},
\end{equation*}
which is a Robin boundary condition for $\phi$. The ratio of these 
two terms is approximately
\begin{equation*}
  \frac{\left|\frac{\phi}{\Delta t^2}\right|}
    {\left|g\frac{\pp{\phi}{n}}{\MM{n}\cdot\MM{k}}\right|}
    \approx \frac{H}{g\Delta t^2} = \left(\frac{H}{c\Delta t}\right)^2,
\end{equation*}
where $c$ is the barotropic wave speed $\sqrt{gH}$. This is the square
of the ratio of the time it takes a barotropic wave to travel a
distance $H$ to the timestep; if we wish to take large timesteps then
this quantity is small and we recover the Neumann boundary condition
\eqref{neumann phi}.

All of this means that if we wish to take large timesteps with an
unstructured mesh, then we must have an efficient method for solving
equation \eqref{phi eqn} with boundary conditions \eqref{neumann phi}.
In this paper we shall see that this equation is ill conditioned when
the aspect ratio $H/L$ is small. We will introduce a new multigrid
preconditioner which allows efficient iterative methods for solving
numerical discretisations on this problem on unstructured meshes.

\subsection{Hydrostatic equations}

In contrast, if one makes the hydrostatic approximation which is used
in many ocean models,
\begin{equation}
  \label{hydrostatic}
  p_z = -g\rho,
\end{equation}
then the pressure can be obtained by integrating this equation, using
the boundary condition $p=0$ on the top surface
$\partial\Omega_{\Top}$. On a vertically structured mesh this equation
can be accurately integrated from top to bottom in columns. On
an unstructured mesh, one approach is to differentiate equation
\eqref{hydrostatic} to obtain the elliptic problem 
\begin{equation}
  \label{hydrostatic elliptic}
  \pp{^2p}{z^2} = -g\pp{\rho}{z}, \quad p=0\,\,\on\,\,\partial\Omega_{\Top},
  \quad \pp{p}{z}=-g\rho\,\,\on\,\,\partial\Omega_{\Floor}.
\end{equation}
The hydrostatic pressure equation \eqref{hydrostatic elliptic} does
not suffer from the same ill conditioning as the nonhydrostatic
pressure equation \eqref{phi eqn} with boundary conditions
\eqref{neumann phi}, since the smallest eigenvalue is independent of
the aspect ratio $\epsilon$ (as we shall see later). This provides the
benchmark for nonhydrostatic pressure equation solver methods, with
the aim of developing methods for the nonhydrostatic equation which
are as fast as methods for the hydrostatic equation.

\section{Finite element formulation}
\label{finite element}

In this paper we consider the finite element approximation to the Poisson
equation,
\begin{equation*}
  -\nabla^2 \phi = f,
\end{equation*}
on the domain $\Omega$ with homogeneous Neumann boundary conditions
\begin{equation*}
  \pp{\phi}{n}(\MM{x}) = 0, \quad \forall \MM{x} \in \Omega.
\end{equation*}
We consider a domain $\Omega$ with horizontal scale $L$, vertical
scale $H$, horizontal velocity scale $U$, and choose nondimensional
domain coordinates
\begin{equation*}
  x' = x/L, \quad y'=y/L, \quad z'=z/H, 
\end{equation*}
In these coordinates, the Poisson equation becomes
\begin{equation}
  \label{rescaled Poisson}-
  \left(\epsilon^2
  \left(\pp{^2}{{x'}^2} + \pp{^2}{{y'}^2}\right)+
  \pp{^2}{{z'}^2}\right)\phi = f', \qquad f' = H^2f,
\end{equation}
with boundary conditions
\begin{equation}
  \label{rescaled bcs}
  (\epsilon^2n'_1,\epsilon^2n'_2,n'_3)\cdot\nabla'\phi = 0, \quad
  \on\,\,\partial\Omega', \quad \MM{n'}=(n'_1,n'_2,n'_3)
\end{equation}
where $\epsilon = H/L$ is the aspect ratio of the unrescaled domain
$\Omega$, and $\Omega'$ is the rescaled domain with normal
$\MM{n}'$. Henceforth we drop all the primes.

Given a solution $\phi$ of equation \eqref{rescaled Poisson} subject
to boundary conditions \eqref{rescaled bcs}, it is possible to obtain
a whole family of solutions
\begin{equation*}
  \phi_c = \phi + c, \quad c \in \mathbb{R},
\end{equation*}
and hence the solution is not unique. However, the constant $c$ does
not have a physical effect since the pressure only appears in the 
Navier Stokes equations as a gradient. This means that we may arbitrarily
fix this constant. This is usually done by requiring that
\begin{equation}
  \label{integral condition}
  \int_\Omega \phi \diff{V} = 0,
\end{equation}
or 
\begin{equation}
  \label{point condition}
  \phi(\MM{x}_0) = 0,
\end{equation}
for some chosen point $\MM{x}_0\in \Omega$. Condition \eqref{point
  condition} is easier to implement in the finite element method on a
general unstructured mesh, and condition \eqref{integral condition}
can subsequently be imposed by subtracting off the integral of $\phi$.
Hence, we shall require condition \eqref{point condition}. Here we
shall additionally require that $\MM{x}_0\in\partial\Omega_{\Top}$;
the approach we describe does not depend on this requirement but it
simplifies the exposition. It was noted in \citet{BoLe05} that good
performance can also be obtained with the Conjugate Gradient method
if the value of the constant mode is not fixed in the assembled 
equations, but instead projected out each 
iteration of the CG solver, but we do not discuss 
that case in this paper.

A finite element discretisation of this equation is obtained by
writing the weak form of equations \eqref{rescaled Poisson} with
boundary conditions \eqref{rescaled bcs}. First we define the 
function space
\begin{equation*}
  H^1(\MM{x}_0) = \{\phi\in H^1: \phi(\MM{x}_0) = 0\},
\end{equation*}
where $H^1$ is the Sobolev space with norm
\begin{equation*}
  \| \phi \|_{H^1}^2 = \int_{\Omega} |\phi|^2 + |\nabla\phi|^2\diff{V},
\end{equation*}
and we seek $\phi\in H^1(\MM{x}_0)$ such that
\begin{equation}
  \label{weak form}
  B_\epsilon(\psi,\phi) = F(\psi), 
\end{equation}
for all test functions $\phi \in H^1(\MM{x}_0)$, where
\begin{eqnarray*}
B_\epsilon(\psi,\phi) &= &\int_\Omega \pp{\psi}{z}\pp{\phi}{z} + \epsilon^2
\left(\pp{\psi}{x}\pp{\phi}{x} + \pp{\psi}{y}\pp{\phi}{y}\right) 
\diff{V} \\
F(\psi) &= & \int_\Omega \psi f \diff{V}.
\end{eqnarray*}

To construct the Galerkin finite element discretisation of these
equations (see \citet{BrSc04} for example), we select a finite
dimensional trial space $V(\MM{x}_0) \subset H^1(\MM{x}_0)$ (typically
by constructing a polygonal mesh on the domain $\Omega$ and
constructing piecewise polynomials), and seek $\phi^\delta\in
V(\MM{x}_0)$ such that
\begin{equation}
  \label{galerkin}
  B_\epsilon(\psi^\delta,\phi^\delta) 
  = F(\psi^\delta),
\end{equation}
for all test functions $\psi^\delta \in V$. If we instead wish to 
solve the equations with a Dirichlet boundary condition imposed on the
top surface $\partial\Omega_{\Top}$, 
\begin{equation}
  \label{dirichlet}
  \phi^\delta = g^\delta, \quad \forall \MM{x}\in \partial\Omega_{\Top},
\end{equation}
then we construct the homogeneous Dirichlet space $\overline{V}$ out
of functions in $V$ which are zero, 
\begin{equation*}
  V(\partial\Omega_{\Top}) = \left\{ \phi^\delta: \phi^\delta \in V,\quad
    \phi^\delta(\MM{x})=0 \quad \forall \MM{x} \in \partial
    \Omega_{\Top} \right\},
\end{equation*}
and seek $\overline{\phi}^\delta\in V(\partial\Omega_{\Top})$ such that
\begin{equation}
  \label{galerkin dirichlet}
  B_\epsilon(\overline{\psi}^\delta,\overline{\phi}^\delta) 
  = F(\overline{\psi}^\delta) - B_\epsilon(\overline{\psi}^\delta,\chi^\delta),
\end{equation}
for all test functions $\overline{\psi}^\delta\in
V(\partial\Omega_{\Top})$, where $\chi^\delta$ is some chosen function
which satisfies \eqref{dirichlet}. The total solution is then
$\phi^\delta = \overline{\phi^\delta} + \chi^\delta$.  This is often
<<<<<<< .mine
referred to as ``lifting'' the boundary conditions; see \citet[for
example]{KaSh2005} for more details. We note that this becomes the
=======
referred to as ``lifting'' the boundary conditions; see 
\citet{KaSh2005} for example for more details. We note that this becomes the
>>>>>>> .r779
Galerkin finite element approximation of the hydrostatic problem
\eqref{hydrostatic elliptic} with suitably chosen $f^\delta$ and when
$\epsilon=0$.

To solve these equations we expand the trial functions $\phi^\delta$,
the test functions $\psi^\delta$, and the right hand side function
$f^\delta$ (or $f^\delta - \chi^\delta$ in the surface Dirichlet case)
in basis function expansions for $V$ (or $\overline{V}$ in the surface
Dirichlet case)
\begin{equation*}
  \phi^\delta(\MM{x}) = \sum_{i=1}^n\phi_i N_i(\MM{x}), \quad
    \psi^\delta(\MM{x}) = \sum_{i=1}^n\psi_i N_i(\MM{x}), \quad
    f^\delta(\MM{x}) = \sum_{i=1}^nf_i N_i(\MM{x}), 
\end{equation*}
and substitute into equation \eqref{galerkin} (or equation
\eqref{galerkin dirichlet}) to obtain a matrix vector equation
\begin{equation}
  \label{Ax=b}
  A_\epsilon\MM{\phi} = M\MM{f},
\end{equation}
where
\begin{equation*}
  A_{\epsilon,ij} = \int_\Omega \pp{N_i}{z}\pp{N_j}{z} + \epsilon^2
    \left(\pp{N_i}{x}\pp{N_j}{x} + 
    \pp{N_i}{y}\pp{N_j}{y}\right) 
    \diff{V}, \quad
    M_{ij} = \int_\Omega N_iN_j
    \diff{V}.
\end{equation*}
Here $A_\epsilon$ is a positive definite sparse matrix, and we wish to
solve equation \eqref{Ax=b} iteratively using the preconditioned
Conjugate Gradient method.  In this paper we shall choose a basis
expansion of $V(\MM{x}_0)$ so that the vector $\MM{\phi}$ of basis
coefficients of $\phi$ takes the form
\begin{equation}
  \MM{\phi} = \begin{pmatrix}
  \MM{\phi}' \\
  \overline{\MM{\phi}} \\
  \end{pmatrix}
  \label{phi_decomposition}
\end{equation}
where $\overline{\MM{\phi}}\in \mathrm{R}^{\overline{m}}$ is the vector of
coefficients corresponding to basis functions which are zero on
$\partial\Omega_{\Top}$, and $\MM{\phi}' \in \mathrm{R}^{m'}$ is the
vector of the coefficients corresponding to the remaining basis
functions. These are typically chosen to vanish on every finite
element ``node'' which is not on $\partial\Omega_{\Top}$.  In this
ordering, we write
\begin{equation}
  A_\epsilon\MM{\phi} = \begin{pmatrix}
  B_\epsilon & C_\epsilon \\
  C^T_\epsilon & \overline{A}_\epsilon \\
  \end{pmatrix}\begin{pmatrix}
  \MM{\phi}' \\
  \overline{\MM{\phi}} \\
  \end{pmatrix},
  \quad M\MM{f} = \begin{pmatrix}
  \MM{b}' \\
  \overline{\MM{b}} \\
  \end{pmatrix}
  \label{equation_decomposition}
\end{equation}
Note that the matrix obtained from equation \eqref{galerkin dirichlet}
(which we denote as $\overline{A}_\epsilon$), is a sub-matrix of the
matrix obtained from equation \eqref{galerkin} (which we denote as
$A_\epsilon$). This matrix solves the problem with Dirichlet boundary
condition at the top surface, rather than Neumann. We shall make
repeated use of this decomposition throughout the rest of the paper.

\section{Condition number estimates}
\label{condition number section}
In this section we obtain estimates in the small aspect-ratio limit
for the condition number of the matrix $A_\epsilon$ which we developed
in the previous section, for both the Neumann and Dirichlet boundary
condition cases. When the condition number is large, the iterative
method can be very slow to converge, and so it is important to
understand the dependence of the condition number $A_\epsilon$ on the
aspect-ratio. We shall note that the Neumann boundary condition case
(which is the case of interest for oceanographic problems) has a
condition number which scales like $\epsilon^{-2}$, whereas the
Dirichlet boundary condition case has a condition number which is
independent of $\epsilon$ as $\epsilon\to 0$. This motivates our
proposed preconditioner.

In this section we estimate the minimum and maximum eigenvalues of
symmetric matrices by using the Rayleigh quotient estimates
\begin{equation*}
  \lambda_{\min} = \min_{\MM{\phi}\neq\MM{0}}\frac{\MM{\phi}^TA_\epsilon\MM{\phi}}
    {\MM{\phi}^T\MM{\phi}}, \qquad
    \lambda_{\max} = \max_{\MM{\phi}\neq\MM{0}}\frac{\MM{\phi}^TA_\epsilon\MM{\phi}}
    {\MM{\phi}^T\MM{\phi}}.
\end{equation*}
To facilitate these estimates, we define the vertical operator $A_0$
and horizontal operator $A_H$ with coefficients
\begin{equation*}
\begin{split}
  A_{0,ij} &= \int_\Omega \pp{N_i}{z}\pp{N_j}{z}
    \diff{V}, \\
  A_{H,ij} &= \int_\Omega
    \pp{N_i}{x}\pp{N_j}{x} + 
    \pp{N_i}{y}\pp{N_j}{y}
    \diff{V},
\end{split}
\end{equation*}
so that
\begin{equation*}
  A_\epsilon = A_0+\epsilon^2 A_H.
\end{equation*}

We first construct an upper bound for the minimum eigenvalue
$\lambda_{\min}$ of $A_\epsilon$. First we note that the intersection of
the null space of $A_0$ with $A_H$ is the zero vector. If $\MM{y}\in
\ker(A_0)$ and $\MM{y} \in \ker(A_H)$, then $\MM{y} \in
\ker(A_\epsilon)$. However, $A_{\epsilon}$ is invertible so
$\MM{y}=0$. This allows us to estimate the minimum eigenvalue:
\begin{eqnarray*}
  \lambda_{\min} &=& \min_{\MM{\phi}\neq\MM{0}}\frac{\MM{\phi}^TA_\epsilon\MM{\phi}}
    {\MM{\phi}^T\MM{\phi}}, \\
  & \leq & \min_{\MM{\phi}\neq\MM{0},\MM{\phi}\in\ker(A_0)}\frac{\MM{\phi}^TA_\epsilon\MM{\phi}}
    {\MM{\phi}^T\MM{\phi}},  \\
  & \leq & \epsilon^2
  \min_{\MM{\phi}\neq\MM{0},\MM{\phi}\in\ker(A_0)}\frac{\MM{\phi}^TA_H\MM{\phi}}
    {\MM{\phi}^T\MM{\phi}}, \\
  & = & c_0\epsilon^2,
\end{eqnarray*}
where 
\begin{equation*}
  c_0 = \min_{\MM{\phi}\neq\MM{0},\MM{\phi}\in\ker(A_0)}\frac{\MM{\phi}^TA_H\MM{\phi}}
    {\MM{\phi}^T\MM{\phi}}, 
\end{equation*}
which is bounded away from zero since $\ker(A_0)\cap
\ker(A_H)=\{\MM{0}\}$.  Next we estimate the maximum eigenvalue
$\lambda_{\max}$ of $A_\epsilon$.
\begin{eqnarray*}
  \lambda_{\max} &=& \max_{\MM{\phi}\neq\MM{0}}\frac{\MM{\phi}^TA_\epsilon\MM{\phi}}
  {\MM{\phi}^T\MM{\phi}}, \\
  & = & \max_{\MM{\phi}\neq\MM{0}}\frac{\MM{\phi}^T\left(A_0
+ \epsilon^2 A_H\right)\MM{\phi}}
  {\MM{\phi}^T\MM{\phi}},  \\
  & \geq & \max_{\MM{\phi}\neq\MM{0}}\frac{\MM{\phi}^TA_0\MM{\phi}}
  {\MM{\phi}^T\MM{\phi}} = c_1,
\end{eqnarray*}
Here, $c_1$ is the maximum eigenvalue of $A_0$ which is independent of
$\epsilon$. Next we compute the condition number of $A_{\epsilon}$
which is the ratio of the largest and smallest eigenvalues
\begin{equation*}
  \Cond(A_\epsilon) = \frac{\lambda_{\max}}{\lambda_{\min}}
  \geq \frac{c_1}{c_0}\epsilon^{-2}.
\end{equation*}
This means that the condition number is unbounded in the small 
aspect ratio limit $\epsilon\to 0$. Since the convergence rate of the
Conjugate Gradient method typically scales with the square root of the
condition number, this means that the Conjugate Gradient method
becomes very slow as $\epsilon \to 0$, and we must find a
preconditioner which makes the condition number independent of
$\epsilon$.

\begin{figure}
\begin{center}
  \includegraphics*[width=10cm]{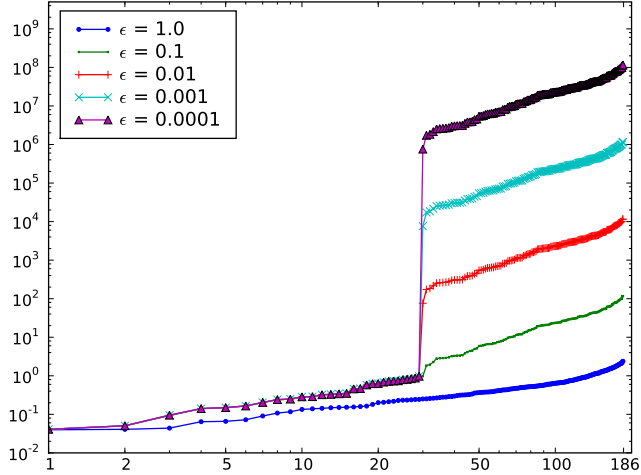} 
  \caption{\label{eigenvalue_lay_hom}Plot showing eigenvalues for the
    matrix $A$ arising from the discretisation of the Poisson equation
    in a box domain with Neumann boundary conditions on all sides. The
    box was decomposed into a tetrahedral mesh divided into a number
    of horizontal layers (\emph{i.e.} the mesh is structured in the
    vertical) and rescaled into various different aspect ratios
    $\epsilon$. Note that there are a cluster of small eigenvalues
    which are independent of $\epsilon$: these eigenvalues correspond
    to the $z$-independent eigenmodes. As $\epsilon$ decreases to
    zero, the width of the spectral gap between these and the
    remaining eigenvalues scales in proportion to $\epsilon^{-2}$.}
\end{center}
\end{figure}

\begin{figure}
\begin{center}
  \includegraphics*[width=10cm]{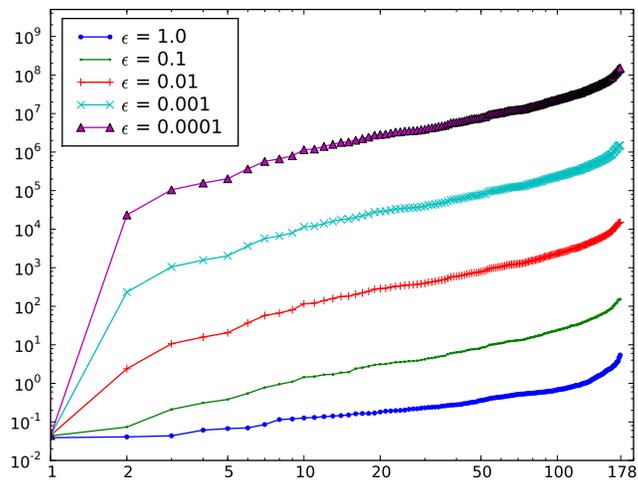} 
  \caption{\label{eigenvalue_lay_uns}Plot showing eigenvalues for the
    matrix $A$ arising from the discretisation of the Poisson equation
    in a box domain with Neumann boundary conditions on all sides. The
    box was decomposed into a tetrahedral mesh which is completely
    unstructured in the vertical and rescaled into various different
    aspect ratios $\epsilon$. Note that there is no longer a
    spectral gap, which is replaced by a more evenly spaced spectrum,
    but that the spread of eigenvalues still increases in proportion
    to $\epsilon^{-2}$.}
\end{center}
\end{figure}

In figures \ref{eigenvalue_lay_hom} and \ref{eigenvalue_lay_uns} we
illustrate these estimates. Two tetrahedral meshes in a box domain
were generated, one which is arranged in horizontal layers, and one
which is fully unstructured in all three dimensions. The finite
element approximation to the Laplace equation was applied to these
meshes, having rescaled the coordinates to various different aspect
ratios. The eigenvalues were then numerically computed using Arnoldi
iteration. Figure \ref{eigenvalue_lay_hom} shows the eigenvalues for
the layered mesh, with various different aspect ratios. We observe a
gap in the spectrum between the eigenvalues corresponding to
$z$-independent eigenvectors (we call these horizontal modes) and the
eigenvalues corresponding to $z$-dependent eigenvectors. The size of
this gap is proportional to $\epsilon^{-2}$. It can also be observed
that the ratio between the largest and smallest eigenvalues is
proportional to $\epsilon^{-2}$. In the unstructured mesh, the
distinction between the horizontal modes and the rest of the
eigenvectors is less clear.  The lack of horizontal alignment in the
mesh means that there are numerical errors in the finite element
approximation of the vertical derivatives which scale with the
horizontal widths $\Delta x$ of the elements (in this case $\Delta
x^2$ since we have used linear finite elements) and hence are the same
order of magnitude (or larger) than the exact eigenvalues of the
horizontal modes. Figure \ref{eigenvalue_lay_hom} shows the
eigenvalues for the unstructured mesh with various different aspect
ratios. We observe the same $\epsilon^{-2}$ scaling for the ratio
between the largest and smallest eigenvalues, but there is no spectral
gap since the small eigenvalues are polluted by numerical errors in
the vertical derivatives.

In contrast, we note very different scaling behaviour when the Neumann
boundary condition on the upper surface (rigid lid) is replaced by a
Dirichlet boundary condition, resulting in the matrix
$\overline{A}_\epsilon$.  For this case, we shall obtain an \emph{upper}
bound on the condition number.  First we bound the minimum eigenvalue
$\overline{\lambda}_{\min}$ of $\overline{A}_\epsilon$ from below.
\begin{eqnarray*}
  \overline{\lambda}_{\min} &=& 
\min_{\overline{\MM{\phi}}\neq\MM{0}}
\frac{\overline{\MM{\phi}}^T\overline{A}_\epsilon\overline{\MM{\phi}}}
  {\overline{\MM{\phi}}^T\overline{\MM{\phi}}}, \\
  & = & \min_{\overline{\MM{\phi}}\neq\MM{0}}
\frac{\overline{\MM{\phi}}^T\left(\overline{A}_0+\epsilon^2\overline{A}_H
    \right)\overline{\MM{\phi}}}
  {\overline{\MM{\phi}}^T\overline{\MM{\phi}}},  \\
  & \geq & \min_{\overline{\MM{\phi}}\neq\MM{0}}
\frac{\overline{\MM{\phi}}^T\overline{A}_0\overline{\MM{\phi}}}
  {\overline{\MM{\phi}}^T\overline{\MM{\phi}}},  \\
  & = & c_2, 
\end{eqnarray*}
which is the minimum eigenvalue of $\overline{A}_0$ and is bounded away
from zero since $\overline{A}_0$ is non-singular.  Next we bound the
maximum eigenvalue $\overline{\lambda}_{\max}$ of $\overline{A}_\epsilon$.
\begin{eqnarray*}
  \overline{\lambda}_{\max} &=& \max_{\overline{\MM{\phi}}\neq\MM{0}}
\frac{\overline{\MM{\phi}}^T\overline{A}_\epsilon\overline{\MM{\phi}}}
  {\overline{\MM{\phi}}^T\overline{\MM{\phi}}}, \\
  & = & \max_{\overline{\MM{\phi}}\neq\MM{0}}
\frac{\overline{\MM{\phi}}^T\left(\overline{A}_0+\epsilon^2\overline{A}_H
\right)\overline{\MM{\phi}}}
  {\overline{\MM{\phi}}^T\overline{\MM{\phi}}},  \\
  & \leq & 2\max_{\overline{\MM{\phi}}\neq\MM{0}}
\frac{\overline{\MM{\phi}}^T\overline{A}_0\overline{\MM{\phi}}}
  {\overline{\MM{\phi}}^T\overline{\MM{\phi}}},  \\
& = & 2c_3, 
\end{eqnarray*}
provided that $\epsilon$ is sufficiently small; here $c_3$ is the
maximum eigenvalue of $\overline{A}_0$. This means that the condition
number of $\overline{A}_\epsilon$ is bounded by
\begin{equation*}
  \Cond(\overline{A}_\epsilon) =
    \frac{\overline{\lambda}_{\max}}{\overline{\lambda}_{\min}} \leq \frac{2c_3}{c_2}
    = 2\Cond(\overline{A}_0),
\end{equation*}
which is twice the condition number of $\overline{A}_0$, and is independent
of $\epsilon$ (this is not a sharp estimate, but illustrates the
scaling with $\epsilon$).

The contrast between the condition number scaling of
$\overline{A}_\epsilon$ and $A_\epsilon$ motivates a preconditioner
strategy in which one solves a reduced problem for the solution on the
surface $\partial\Omega_{\Top}$, and then uses this surface solution
as a Dirichlet boundary condition to reconstruct a solution throughout
$\Omega$. This reconstruction step amounts to inverting
$\overline{A}_\epsilon$ which, as we have just seen, has a condition number
which is independent of $\epsilon$. We shall describe this
preconditioner strategy in the following section.

\section{New preconditioner}
\label{new preconditioner}
In this section we develop our new preconditioner for equation
\eqref{Ax=b}, which is derived from the strategy of eliminating the
degrees of freedom associated with the solution in the interior of
$\Omega$ to give a reduced equation on the surface
$\partial\Omega_{\Top}$.  In each iteration the preconditioner will
approximately solve this problem, and then approximately reconstruct
the solution in the interior using $\overline{A}_\epsilon$. We describe the
preconditioner as follows: in Section \ref{amg section} we briefly
summarise the general algebraic multigrid preconditioning strategy.
In Section \ref{schur section} we introduce a reformulation of
equation the \eqref{Ax=b} that decomposes the inverse of $A$ into a
vertically lumped system and a system with a Dirichlet boundary
condition on top. Sections \ref{extrapolation section} and \ref{sspai
  section} explain the approximations that need to be made to this
decomposition to apply it as a preconditioner.

\subsection{Algebraic multigrid preconditioners}
\label{amg section}
The general idea of multigrid methods is to tackle multiscale, ill
conditioned problems by trying to solve for the different components
of the solution, associated with different length scales, separately.
This is accomplished by a sequence of coarsening operations, in which
the dimension of the problem is reduced step by step. The coarser
system no longer supports the smaller scale features and has therefore
an improved condition number. Thus the large scale, small eigenvalue
modes can be efficiently solved on a reduced system, whereas the small
scale, large eigenvalue modes are easily reduced with standard
preconditioners such as SOR (therefore in this context referred to as
\emph{smoothers}).

Classical geometric multigrid methods, implement this coarsening step
via a coarsening of the mesh on which the problem is defined, for
instance via a $h\to 2h$ coarsening on structured meshes.
\emph{Algebraic multigrid} (AMG) methods (see \citet{stueben01} for an
introduction), use the algebraic properties, matrix graph and
coefficients, of the matrix to construct a coarsening operator. This
more general approach has the advantage that it works equally well for
unstructured mesh discretisations. Additionally it is possible to take
anisotropies in the problem into account by selecting only matrix
graph connections associated with large matrix coefficients.

For symmetric problems, to keep the problem 
symmetric at each level, the prolongation operator $P$,
that maps the solution of the reduced system back to the 
previous level, is usually the transpose of the coarsening 
operator $P^T$. Given the original matrix $A$ at the finest level, 
the matrix of the reduced system, is given by
\begin{equation}
  \label{reduced AMG}
  P^TA P\MM{y} = P^T\MM{b}, \qquad \MM{y} \in \mathbb{R}^m,
\end{equation}
After solving at the coarsest level, the solution is mapped back
using $y\to x=Py$.

The smaller scale modes that are only represented 
in the full problem are reduced by applying a smoother $S$. As this 
smoothing step will also again change the solution at the coarse level, 
\emph{i.e.} the solve at the coarse level and the smoothing for the full problem
are not independent, the whole procedure of restriction, 
coarse solve, prolongation, and smoothing needs to be applied in an 
iterative manner. To keep things symmetric, 
the smoothing step $S$ after prolongation is usually mirrored by a 
transpose smoother $S^T$ before the reduction. For instance a forward 
sweep of SOR before the restriction can be accompanied by a 
backward sweep after the prolongation.

A typical 2-level multigrid
cycle (V-cycle) then looks like
\begin{equation*}
\xymatrix{
  \mathbb{R}^n \ar[r]_{S^T} & \mathbb{R}^n \ar[dr]_{P^T} & & & & \mathbb{R}^n \ar[r]_{S} & \mathbb{R}^n \\
  & & \mathbb{R}^m \ar[rr]_{\quad\quad\left(P^T A P\right)^{-1}\quad\quad} & & \mathbb{R}^m \ar[ru]_{P} \\
}
\end{equation*}
Finally by replacing the coarse solve with a multigrid V-cycle applied 
to the reduced system, the multigrid method can be extended recursively
to multiple levels.

Best results are obtained if the multigrid V-cycle is embedded, 
as a preconditioner, in a Krylov subspace method. Different multigrid 
preconditioning approaches are formed by different coarsening strategies
and different choices of smoothers. The algebraic multigrid preconditioner
used in the results section, implements the smoothed aggregation 
approach of \citet{vanek96}. This method is known to work very well 
for strongly anisotropic elliptic problems. As will be shown in the results
section however the convergence rate is not independent of the 
aspect ratio. Therefore, we seek to improve upon this purely 
algebraic method in the following sections.

\subsection{Schur complement equation}
\label{schur section}
Motivated by the analysis in section \ref{condition number section}, 
where we observed that the condition number of the linear system 
becomes independent of the aspect ratio if the Neumann boundary condition 
on top is replaced by a Dirichlet condition, we proceed by constructing a
reduced system where we first solve for the solution $\MM{\phi}'$ on 
$\partial\Omega_{\text{Top}}$, and then reconstruct
$\overline{\MM{\phi}}$ (cf. the decomposition 
of $\MM{\phi}$ into $\MM{\phi}'$ and $\MM{\overline\phi}$ in
\eqref{phi_decomposition} and \eqref{equation_decomposition}). This can be done by solving the Schur complement
equation
\begin{equation}
  \label{schur}
  \underbrace{(B-C\overline{A}_\epsilon^{-1}C^T)}_{\mbox{Schur matrix}}\MM{\phi}' = \MM{b}' 
    - C\overline{A}_\epsilon^{-1}\overline{\MM{b}},
\end{equation}
and then solving
\begin{equation}
  \label{reconstruction}
  \overline{A}_\epsilon\overline{\MM{\phi}} = -C^T\MM{\phi}' + \MM{b}',
\end{equation}
to reconstruct the solution in the interior. Note that
$\dim(\MM{\phi}') \ll \dim(\MM{\phi})$, and also that the
reconstruction equation \eqref{reconstruction} has a condition number
which is independent of $\epsilon$ as $\epsilon\to 0$.

The Schur complement matrix contains $\overline{A}_\epsilon^{-1}$ and hence
is a full matrix which is expensive to assemble and solve. Hence we
shall propose a strategy to form approximations to equations
\eqref{schur} and \eqref{reconstruction} which can be used as a
preconditioner in a manner similar to a multigrid preconditioner.

\subsection{Extrapolation operator}
\label{extrapolation section}
To build the approximation to the Schur matrix, we first note that
equation \ref{schur} may be rewritten as
\begin{equation}
  \label{E equation}
  E^TA_{\epsilon}E\MM{\phi}' = E^T\MM{b},
\end{equation}
where
\begin{equation*}
  E = 
  \begin{pmatrix}
    I \\
    -\overline{A}_\epsilon^{-1}C^T \\
  \end{pmatrix}.
\end{equation*}
We call $E$ the extrapolation operator. Given $\MM{\phi}'\in \mathbb{R}^m$,
the operation
\begin{equation*}
  \MM{\phi} = E\MM{\phi}',
\end{equation*}
produces the finite element discretisation of the solution of 
the Laplace equation
\begin{equation}
  \label{laplace}
  \nabla^2\phi = 0,
\end{equation}
with Neumann boundary conditions $\partial\psi/\partial n=0$ on the
coasts and bottom surface, and Dirichlet boundary conditions
\begin{equation}
  \label{dirichlet_top}
  \phi = \phi'
\end{equation}
on the top surface, where $\phi$ is the function on $\Omega$ with
finite element basis function coefficient vector $\MM{\phi}$ and
$\phi'$ is the function on $\partial\Omega_{\Top}$ with finite element
basis function coefficient vector $\MM{\phi}'$. Note that in the
small aspect ratio limit, \eqnref{laplace} converges to
\begin{equation*}
  \pp{^2}{z^2}\phi = 0,
\end{equation*}
with $\phi=\phi'$ on the top surface, and $\partial\phi/\partial z=0$
on the bottom surface. The solution $\phi$ is then the vertical
extrapolation of $\phi'$ \emph{i.e.} 
\begin{equation*}
  \phi(x,y,z) = \phi'(x,y).
\end{equation*}
We shall use this in subsequent sections to construct an approximation
to $E$.

We can eliminate $\overline{\MM{\phi}}$ using equation
\eqref{reconstruction} to obtain
\begin{eqnarray*}
\MM{\phi} &=& \begin{pmatrix}
\MM{\phi}' \\
\overline{\MM{\phi}} \\
\end{pmatrix} \\
& = & 
\begin{pmatrix}
\MM{\phi}'
\\
\overline{A}_\epsilon^{-1}\left(-C^T\MM{\phi}'+\MM{b}'\right)
\end{pmatrix} \\
& = & 
E\MM{\phi}' +
\begin{pmatrix}
0 \\
\overline{A}_\epsilon^{-1}\MM{b}' \\
\end{pmatrix} \\
& = & E(E^TA_\epsilon E)^{-1}E^T\MM{b} + 
\begin{pmatrix}
\MM{0} \\
I \\
\end{pmatrix}  
\overline{A}_\epsilon^{-1}
\begin{pmatrix}
\MM{0} \\
I \\
\end{pmatrix}^T\MM{b},
\end{eqnarray*}
that is
\begin{equation}
  \label{E multigrid}
  A_\epsilon^{-1}=
    \left(
    E(E^TA_\epsilon E)^{-1}E^T + 
    \underbrace{
    \begin{pmatrix}
      \MM{0} \\
      I \\
    \end{pmatrix}  
    \overline{A}_\epsilon^{-1}
    \begin{pmatrix}
      \MM{0} \\
      I \\
    \end{pmatrix}^T}_{\mathrm{smoother}}
    \right).
\end{equation}
The two parts of this formula 
may be interpreted in the context of a general
multigrid strategy. The first term is similar to a 2-level 
multigrid cycle with prolongation operator $E$. It projects the equation 
to a vertically lumped system. The second term acts on the vertical modes 
in the solution and can therefore be seen as an additive smoother. It is 
to be noted that \eqref{E multigrid} provides an exact solution for 
the inverse of $A_\epsilon$, provided both inversions are performed 
exactly. However, both $E$ and $\overline{A}_\epsilon^{-1}$ are
dense matrices; in the next two sections we will provide approximations for
both the extrapolation operator $E$ and the inverse of $\overline{A}_\epsilon$, 
such that \eqref{E multigrid} can be used as a preconditioner for equation 
\eqref{Ax=b}.

\begin{figure}[htb]
\begin{center}
\includegraphics[width=0.7\textwidth]{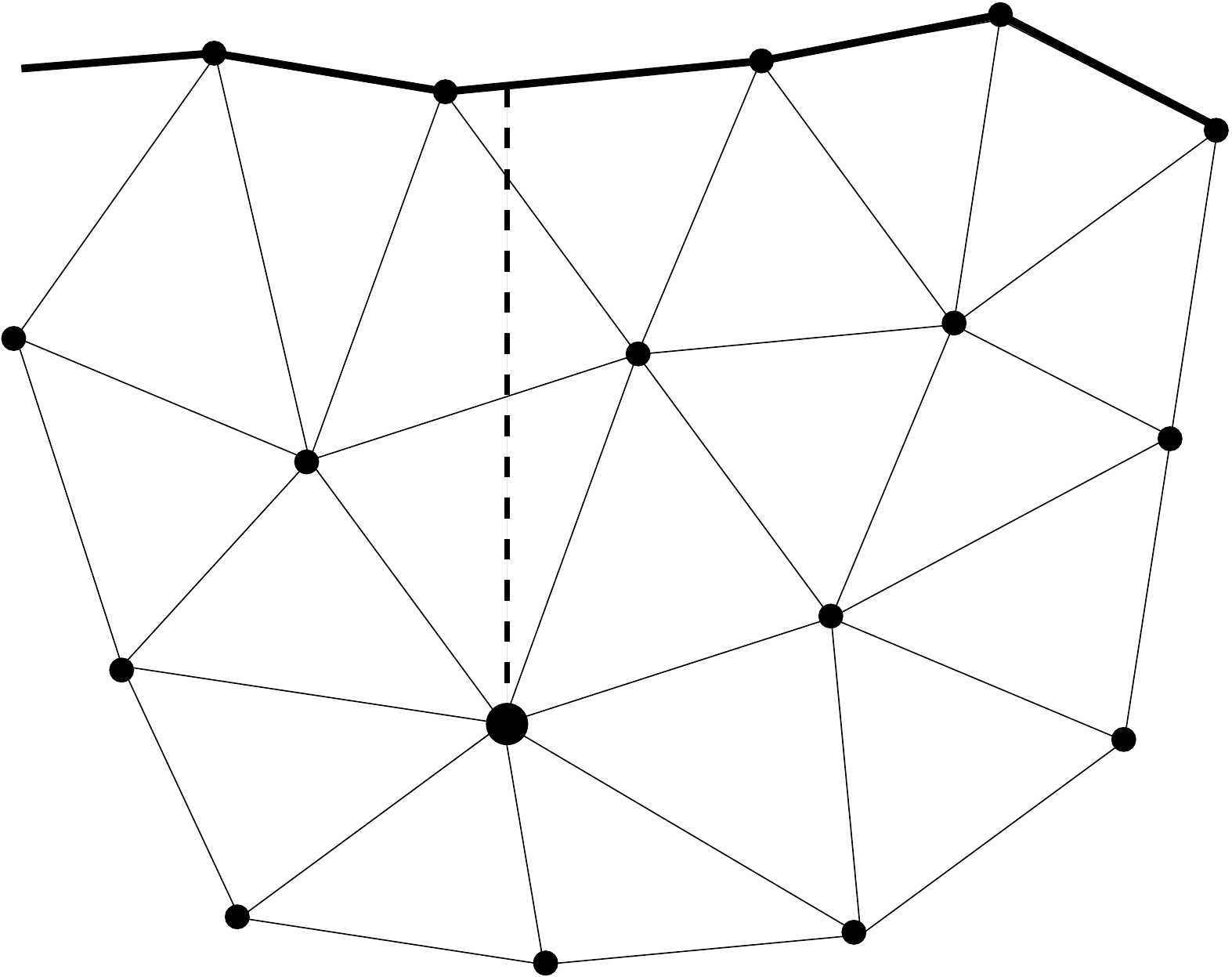}
\end{center}
\caption{Vertical extrapolation performed by projecting each node in 
the full mesh (here depicted in 2D) straight upward onto the 
top surface mesh (here 1D), and interpolating the value in this 
projected node from the surrounding nodes in the surface element.}
\label{fig:vertical_interpolation}
\end{figure}
\subsection{Approximation of the vertical extrapolation operator $E$}
\label{sspai section}
As noted in the previous section, the operator $E$ given by
discretisation of \eqref{laplace} with Dirichlet boundary conditions
condition \eqref{dirichlet_top} on top, converges to a vertical 
extrapolation operator as $\epsilon\to 0$. We therefore expect that 
for large $\epsilon$, such a vertical extrapolation operator, 
between the top surface mesh and the full mesh, 
is a good approximation of $E$. This operator can simply be 
constructed by projecting nodes of the full mesh
in the vertical direction onto the surface mesh, and interpolating 
within the surface triangle each projected node lies 
(see figure \ref{fig:vertical_interpolation}).
This gives an 
approximation $\tilde E:\mathrm{R}^{m'}\to\mathrm{R}^{\overline{m}}$ of $E$
with a limited stencil: for a continuous linear (P1) discretisation in 
3 dimensions it connects every interior node with three nodes of the 
surface mesh.

Another approach to find a sparse approximation of $E$ is given 
by projecting $E$ onto a chosen sparsity pattern using a 
modification of the symmetric sparse
approximate inverse (SSPAI) \citep{BeMeTu1996}. $\tilde{E}$ is obtained
by writing
\begin{equation*}
  \tilde{E} = \begin{pmatrix}
    I \\
    F \\
  \end{pmatrix}
\end{equation*}
and then selecting a sparsity pattern for $F$. The non sparse entries
of $F$ are then obtained by minimising
\begin{equation*}
  \|F^T\overline{A}_\epsilon - C\|^2,
\end{equation*}
subject to the sparsity constraints, where $\|\cdot\|$ is the
Frobenius norm. This leads to decoupled sparse matrix problems
\begin{equation*}
  {\overline{A}}_i\MM{v}_i = -\MM{r}_i, \qquad i=1,\ldots, n,
\end{equation*}
where $\MM{v}_i$ is the $i$-th column of $F$ restricted to the
sparsity pattern of that column, $\MM{r}_i$ is the $i$-th column of
$C^T$ restricted to the sparsity pattern, and $\overline{A}_i$ is
$\overline{A}_\epsilon$ restricted to the sparsity pattern of the $i$-th
column. We do not pursue this approach in this paper, preferring to
use the projection approach described above.

\subsection{Additive smoother}
\label{smoother section}

By replacing $E$ with $\tilde{E}$, we have produced the approximate
inverse
\begin{equation}
  A_\epsilon^{-1} \approx \tilde{E}(\tilde{E}^TA_\epsilon\tilde{E})^{-1}\tilde{E}^T 
    + 
    \begin{pmatrix}
      \MM{0} \\
      I \\
    \end{pmatrix}
    \overline{A}_\epsilon^{-1}
    \begin{pmatrix}
      \MM{0} \\
      I \\
    \end{pmatrix}^T.
  \label{ourpreconditioner}
\end{equation}
To use this as a preconditioner, we must also approximate the additive
smoother. The second term could be evaluated exactly by solving a
matrix equation $\overline{A}_\epsilon\overline\phi=\overline b$ for the interior
part of the residual. Although, as noted before, this system is much
better conditioned than the full system, the solution of an elliptic
equation on the interior of the mesh is still quite an expensive
operation that needs to be performed during each application of the
preconditioner within each the Krylov iteration. Moreover the solution
of this interior equation needs to be done using an iterative Krylov
method as well. It is well known that embedding a Krylov method within
another Krylov iteration, requires the use of a flexible Krylov method
for the outer iteration (\emph{e.g.} FGMRES \citep{saad93}). A major drawback would
therefore be that this approach would inhibit the use of the Conjugate
Gradient method for the outer iteration.

A very simple smoothing strategy is obtained by realising that the
first term of the proposed preconditioner is just the first stage of a
general multigrid method. In this projection the long scale,
horizontal modes are separated out and the vertical projection can
therefore be interpreted as a general coarsening step such as those in
any multigrid method. The necessary smoothing step to filter out the
short scale modes, is there often done by application of one or more
SOR iterations of the entire system. For small aspect ratio problems
this may therefore be enough to reduce the vertical modes in the
error.

In some cases, the mesh may contain a lot of structure in the vertical
as well as in the horizontal. For instance an adaptive mesh model
might focus resolution on physics related to the baroclinic modes of
the system. In such cases the simple SOR smoothing may not be
enough. The vertical lumping step would take out too much structure in
one step. This may be compared to so called ``aggressive coarsening''
techniques in general multigrid methods that are usually accompanied
with improved smoothing techniques. A more accurate approximation of
the second term in \eqref{ourpreconditioner} would be to replace the
inverse matrix $\overline{A}_\epsilon^{-1}$ by a full cycle of the general
AMG method applied to $\overline{A}_\epsilon$. The next section will
provide a comparison of the simple SOR smoother with this more
advanced additive smoother.

\section{Numerical experiments}

\label{numerics}

In this section we present numerical results which test out our
preconditioner on matrices obtained from the linear finite element
approximation of the Laplace equation with the horizontal coordinates
rescaled to various aspect ratios. The solvers were developed using
the open source PETSc library \citep{petsc-efficient}.

To compute errors, we selected a right hand side for the matrix vector
equation by choosing a solution and multiplying it by the matrix. This
allows us to compute errors exactly at each iteration of the
solvers. Throughout this section we use the $\inf$-norm to measure the
magnitude of the error: our rationale for this is that we are
motivated by multiscale applications in which one may be very
concerned with the numerical solution in one small region of the
domain (for example one may wish to embed a convection cell in an
ocean basin and observe how it is affected by the large scale
dynamics). In this case it may be possible to obtain a small $L_2$
error whilst the solution in the small region is still inaccurate. We
present plots of error against number of iterations, and also against
floating point operations (flops). The flop count is provided by a
intrinsic PETSc routine.

In this section we obtain results from two meshes, both of a
1$\times$1$\times$1 cube; the coordinates of the meshes are then
rescaled to a range of aspect ratios in the small aspect ratio
limit. Mesh A is a Delaunay triangulation for a set of 57453 roughly
equispaced points in the cube; the points are not arranged in layers
and hence the mesh is unstructured in the vertical. Mesh B is a
Delaunay triangulation on a mesh in which the majority of the points
are clustered at the centre of the cube (there are 99017 points in
total), leading to very small elements. This is a truly multiscale
mesh in which small eigenvalues exist due to both the small aspect
ratio and also due to small elements. Our aim is to develop robust,
efficient matrix solvers for these challenging multiscale meshes.

The three preconditioners that are compared are:
\begin{itemize}
\item A general AMG method based on the \emph{smoothed aggregation
    method} \citep{vanek96}. This uses our own
  implementation constructed using the ``MG'' interface provided by
  PETSc. The smoothing at each level is given by a single forward SOR
  sweep ($\omega=1.0$) as a pre-smoother and a backward sweep for
  post-smoothing.  The coarsening strategy is based on the
  strongly-coupled connection criterion

\begin{equation*}
  |A_{ij}|>\varepsilon\sqrt{A_{ii} A_{jj}}
\end{equation*}

where a $\varepsilon$ of $0.01$ has been chosen. The smoothing in 
the aggregation operator uses $\omega=2/3$.
\item
The preconditioner given by the \emph{vertically lumped approach}

\begin{equation*}
  A_\epsilon^{-1} \approx 
    \tilde{E}(\tilde{E}^TA_\epsilon\tilde{E})^{-1}\tilde{E}^T
\end{equation*}

where the vertically lumped system is approximately
solved using a single multigrid cycle applied 
to $\tilde{E}^TA_\epsilon\tilde{E}$. This is combined with a single 
forward and backward SOR sweep as respectively a pre and post 
smoothing step. Thus the vertical lumping operator $\tilde E$ is 
treated as an ordinary coarsening operator, and the vertical lumping 
of the equation is simply the first of a multilevel multigrid 
cycle.

\item
As a last approach, the above multigrid cycle, 
including the vertical lumping as the first coarsening step, 
is combined with the \emph{additive smoother}, where $\overline{A_\epsilon}^{-1}$ 
is approximated applying a cycle of the smoothed aggregation AMG method 
to $\overline{A_\epsilon}$.
\end{itemize}
In all cases the full multigrid cycle is applied as a preconditioner 
within each iteration of the Conjugate Gradient method.

\subsection{General multigrid methods}
\begin{figure}
\begin{center}
\includegraphics*[width=12cm]{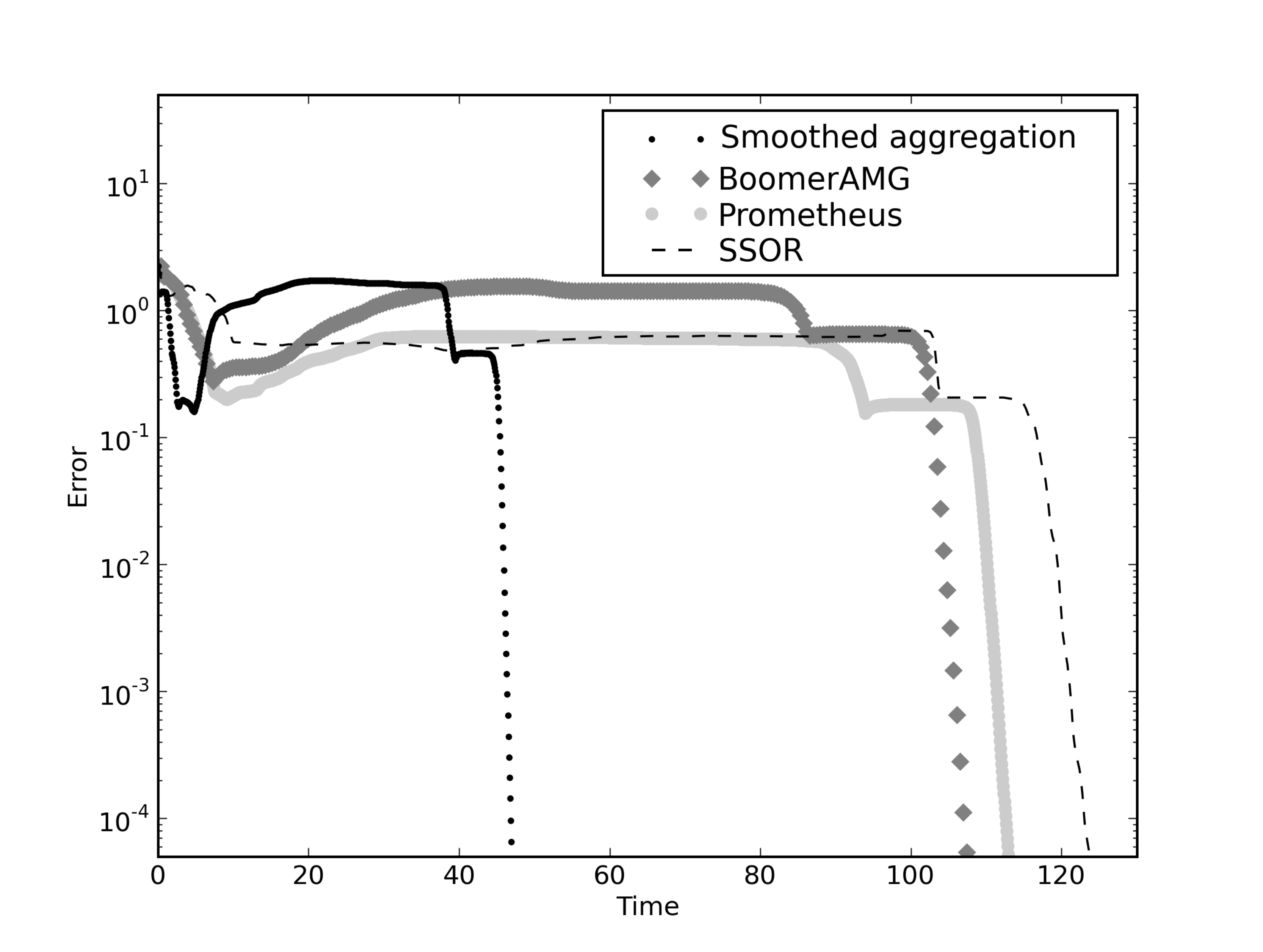}
\end{center}
\caption{\label{other_methods_mesha} Plot showing error versus
  computation time for the smoothed aggregation multigrid
  preconditioner combined with the Conjugate Gradient method, compared
  with two other multigrid preconditioners available within PETSc
  (namely BoomerAMG and Prometheus), and the SSOR preconditioner. The
  smoothed aggregation method uses a coarsening strategy which is
  weighted by the matrix entries. When the aspect ratio is small (it
  is $1/1000$ in this case), this method tends to aggregate degrees of
  freedom which are nearly vertically aligned. The smoothed
  aggregation method is substantially more effective than the other
  multigrid methods and the SSOR method, but we still observe a long
  period during which the error is not reduced.  }
\end{figure}

In figure \ref{other_methods_mesha}, we compare the smoothed
aggregation preconditioner to two other algebraic multigrid algorithms
available through the PETSc library (BoomerAMG and Prometheus), and
the classical Symmetric Successive Over-Relaxation (SSOR)
preconditioner. The four preconditioners are combined with the
Conjugate Gradient method in solving the pressure Poisson equation on
Mesh A rescaled to an aspect ratio of $1/1000$ (a reasonable aspect ratio
for large scale oceanographic applications).  It is shown that the
smoothed aggregation is substantially more effective than the other
methods (we additionally notice that the other multigrid methods are
not much more effective than SSOR in this small aspect ratio example,
although both methods contain a number of parameters and it may be
possible to obtain better results by tuning). All four methods produce
long ``plateaus'' in the error that are maintained for many iterations
before the error finally drops. This means that none of these
preconditioners result in feasible methods for solving the pressure
Poisson equation in small aspect ratio domains.

\begin{figure}
\begin{center}
\includegraphics*[width=12cm]{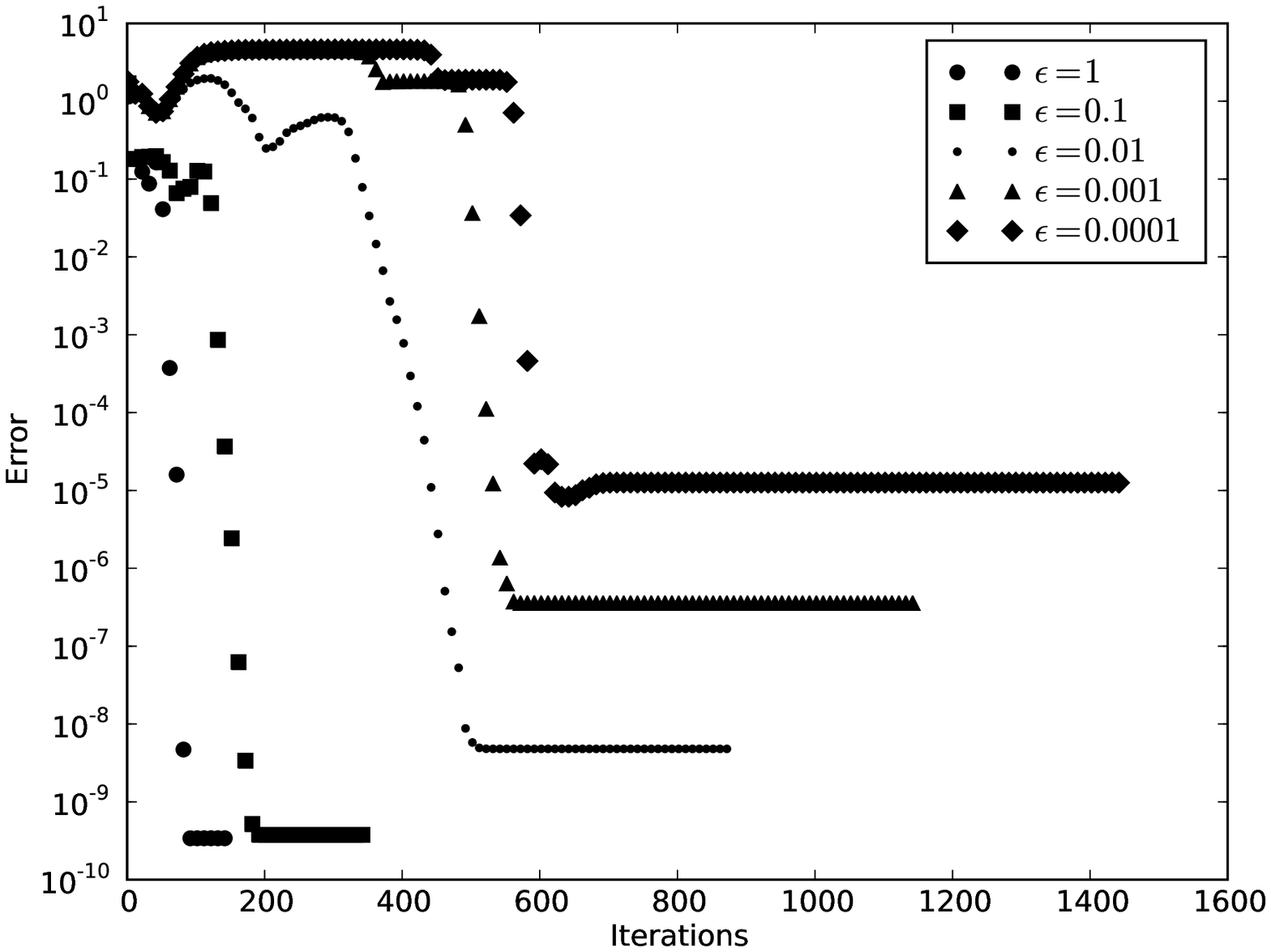}
\end{center}
\caption{\label{iso_no_vl_eps} Plot showing error (using the
  $\infty$-norm) against number of iterations, for the Conjugate
  Gradient method applied to the Poisson equation discretised on mesh
  A, using the smoothed aggregation multigrid preconditioner. The mesh has been
  rescaled to various different aspect ratios $\epsilon$ as
  indicated in the plot.  The number of iterations required to
  converge increases with decreasing $\epsilon$, with a long ``plateau'' for
  small aspect ratios.}
\end{figure}

In figure \ref{iso_no_vl_eps}, we further illustrate the problems that
occur when the Conjugate Gradient method with the smoothed aggregation
multigrid preconditioner is applied to the Poisson equation on mesh
A. The error is plotted against the iteration number for the
preconditioned Conjugate Gradient method for various aspect ratios. As
the aspect ratio of the domain decreases, the condition number
increases and the convergence rate of the iterative solver gets
slower. We observe a ``plateau'' in the error which is maintained for
many iterations before the error finally drops; this plateau becomes
longer and longer as the aspect ratio decreases. For small aspect
ratios this makes the iterative solver prohibitively slow.

\subsection{Preconditioner with vertical lumping}
\begin{figure}
\begin{center}
\includegraphics*[width=12cm]{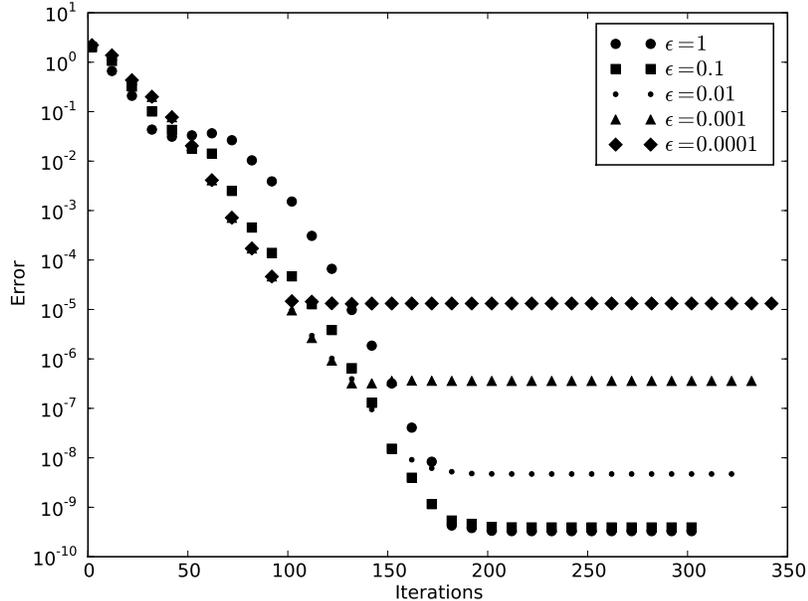}
\end{center}
\caption{\label{iso_vl_eps} Plot showing error (using the
  $\infty$-norm) against number of iterations, for the Conjugate
  Gradient method applied to the Poisson equation discretised on mesh
  A, using the vertically lumped preconditioner. The mesh has been
  rescaled to various different aspect ratios $\epsilon$ as
  indicated in the plot.  The convergence rate becomes independent of
  $\epsilon$ for small aspect ratios.}
\end{figure}

In figure \ref{iso_vl_eps}, the same information (error plotted
against iteration number for the solution of the matrix system
obtained from Mesh A) is given for the vertically lumped
preconditioner using an SOR smoother. We note that in contrast to the
standard multigrid preconditioners tested in the previous subsection,
there is no plateau and the convergence rate becomes independent of
$\epsilon$ for small aspect ratios. We attribute this fast convergence
to the removal of small eigenvalues in (nearly) vertically-independent
eigenmodes by the vertically lumped preconditioner. For small aspect
ratios there is an exponential decay of error with iteration from the
very first iterations. 

As an aside, we observe that the remaining error in the approximation
after the solver has stopped converging, increases with decreasing
$\epsilon$. We ascribe this to numerical round off error (all runs are 
done in double precision). The scaling 
of the condition number with $\epsilon$ is consistent with the observed
loss in accuracy. The smallest used $\epsilon$ of $0.0001$ still gives 
an accuracy that is acceptable. This remaining error will show 
up in all further figures.

\begin{figure}
\begin{center}
\includegraphics*[width=12cm]{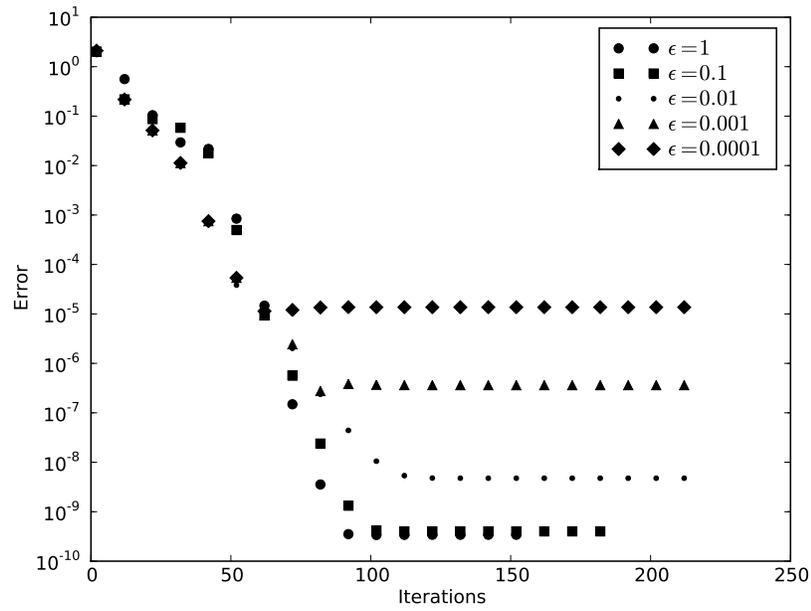}
\end{center}
\caption{\label{iso_vl_as_eps} Plot showing error (using the
  $\infty$-norm) against number of iterations, for the Conjugate
  Gradient method applied to the Poisson equation discretised on mesh
  A, using the vertically lumped preconditioner combined with additive
  smoothing. The mesh has been rescaled to various different
  aspect ratios $\epsilon$ as indicated in the plot.  The
  convergence rate becomes independent of $\epsilon$ for small
  aspect ratios, but the smoother does not improve the convergence
  much for this mesh.}
\end{figure}

In figure \ref{iso_vl_as_eps}, the error is plotted against iteration
for the vertically lumped preconditioner using our additive
smoother. We note that the error again decays exponentially with
iteration number at a rate which is independent of $\epsilon$ for
small aspect ratios. However, one sweep of the additive smoother is
more expensive than one SOR sweep, and hence it is necessary to
compare the performance of the two smoothing strategies in terms of 
computational cost as well as number of iterations.

\begin{figure}
\begin{center}
\includegraphics*[width=12cm]{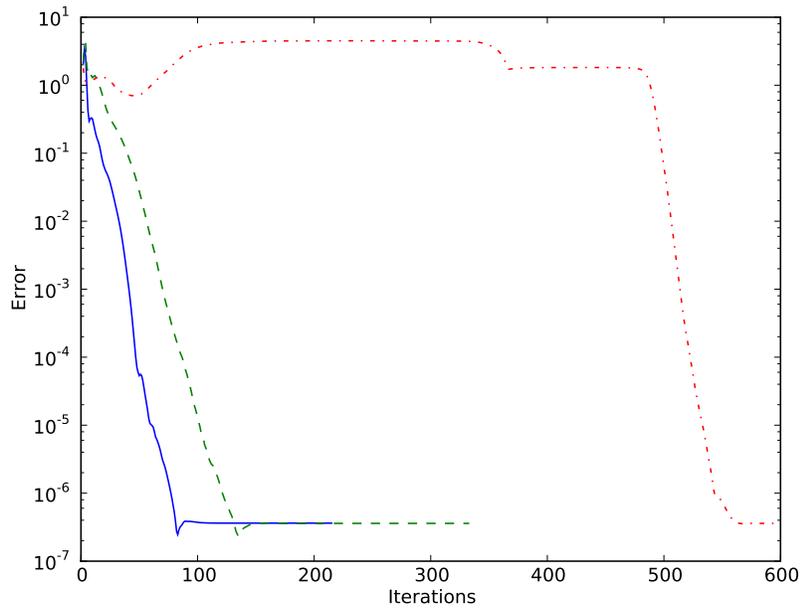}
\end{center}
\caption{\label{iso_eps0_001_its} Plot showing error (using the
  $\infty$-norm) against number of iterations, for the Conjugate
  Gradient method applied to the Poisson equation discretised on mesh
  A ($\epsilon=0.001$), with various different preconditioners. The
  continuous line indicates the vertically lumped preconditioner with
  the additive smoother, the dashed line indicates the vertically
  lumped preconditioner without the additive smoother, and the
  dash-dotted line indicates the smoothed aggregation multigrid preconditioner.}
\end{figure}

\begin{figure}
\begin{center}
\includegraphics*[width=12cm]{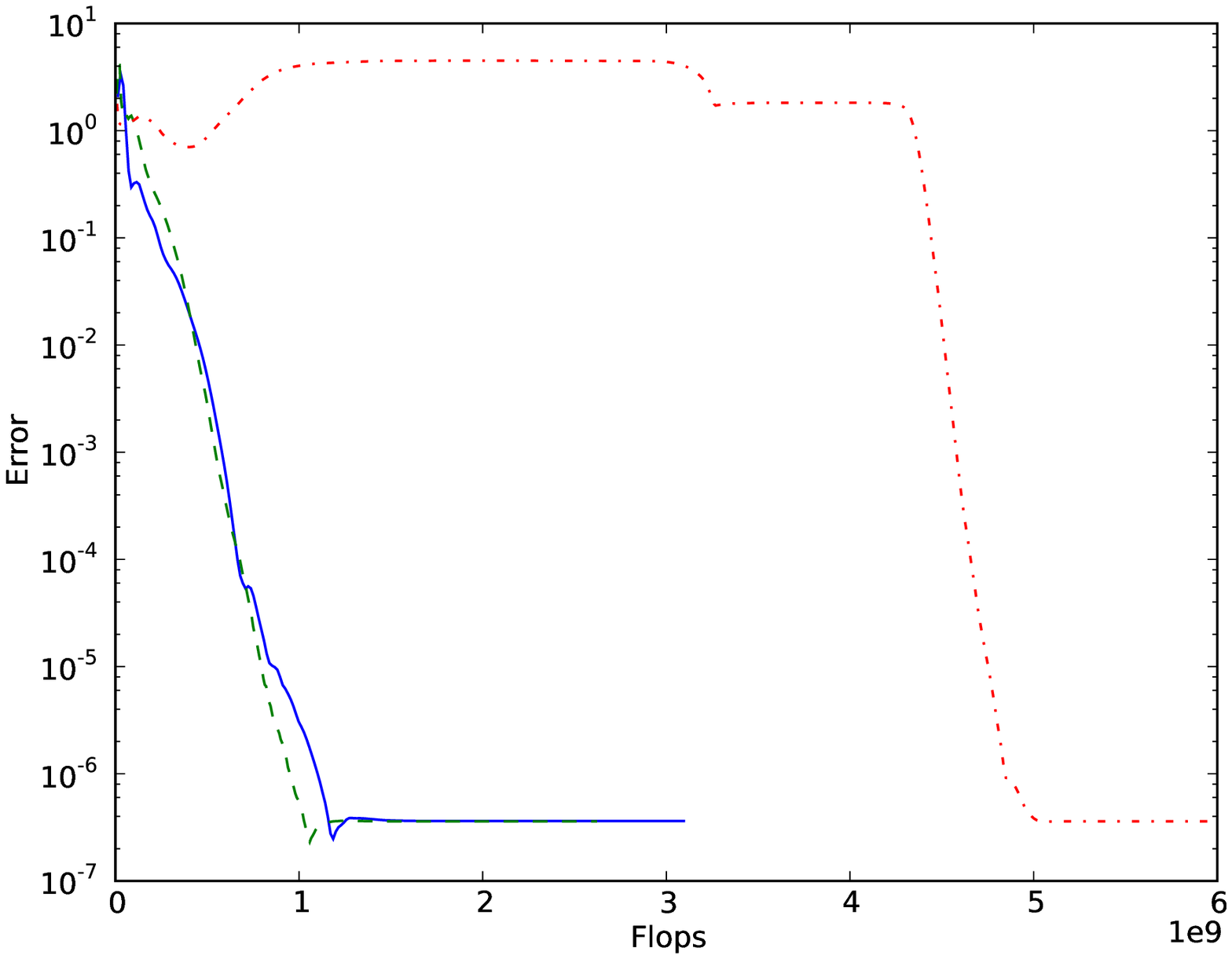}
\end{center}
\caption{\label{iso_eps0_001_flops} Plot showing error (using the
  $\infty$-norm) against number of floating point operations (flops)
  counted by PETSc, for the Conjugate Gradient method applied to the
  Poisson equation discretised on Mesh A ($\epsilon=0.001$), with various different
  preconditioners. The continuous line indicates the vertically lumped
  preconditioner with the additive smoother, the dashed line indicates
  the vertically lumped preconditioner without the additive smoother,
  and the dash-dotted line indicates the smoothed aggregation multigrid
  preconditioner.}
\end{figure}

In figure \ref{iso_eps0_001_its}, the error is plotted against number
of iterations for the smoothed aggregation preconditioner and the
vertically lumped preconditioner with and without the additive
smoother for the matrix obtained from Mesh A with $\epsilon=0.001$.
We observe that the vertically lumped preconditioner converges much
faster than the smoothed aggregation preconditioner, and that the
additive smoother reduces the number of iterations required for
convergence. The vertically lumped preconditioner has made it feasible
to solve the pressure Poisson equation on this type of mesh.  However,
as we observe in figure \ref{iso_eps0_001_flops} which shows the error
plotted against flops, the additive smoother is more expensive than
SOR since it involves several applications of SOR at different
levels. We observe that for Mesh A, the time to converge is
approximately the same with the SOR smoother or with the additive
smoother. The vertically lumped preconditioner produces an
approximation to the vertically-independent (barotropic) component of
the solution with very small eigenvalues and it is the job of the
smoothers to approximate the vertically-varying (baroclinic)
component. These results show that for Mesh A, which has roughly
isotropic tetrahedra before rescaling to small aspect ratio, the SOR
smoother is reasonably effective in approximate the baroclinic
components.

\begin{figure}
\begin{center}
\includegraphics*[width=10cm]{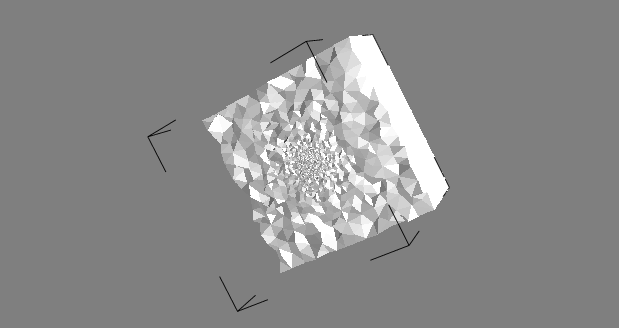}
\end{center}
\caption{\label{monster} Plot showing ``cutaway'' surface through mesh
  B, used for benchmarking the
  preconditioner. The mesh has very fine mesh elements at the
  middle of the domain, so that there are a large range of
  lengthscales in the mesh.}
\end{figure}

\begin{figure}
\begin{center}
\includegraphics*[width=12cm]{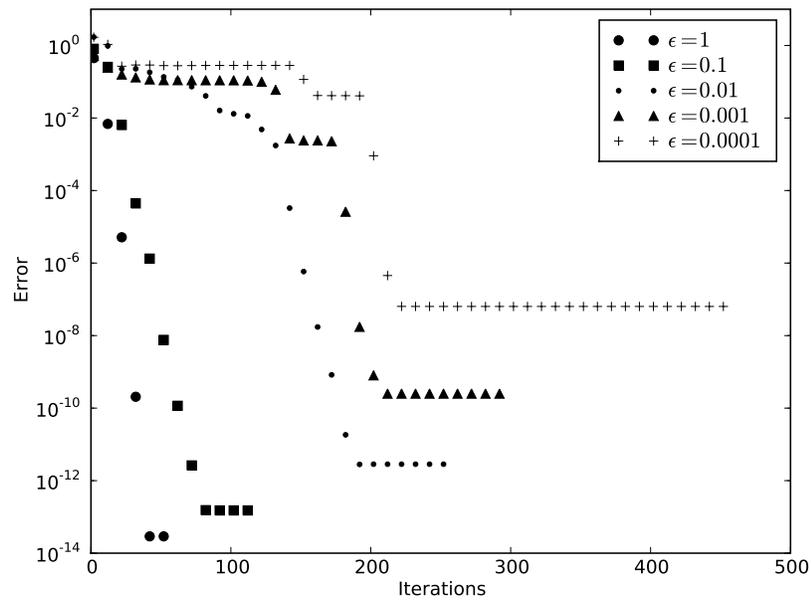}
\end{center}
\caption{\label{monster_no_vl_eps} Plot showing error (using the
  $\infty$-norm) against number of iterations, for the Conjugate
  Gradient method applied to the Poisson equation discretised on mesh
  B (shown in figure \ref{monster}, $\epsilon=0.001$), using the smoothed aggregation multigrid
  preconditioner. The mesh has been rescaled to various different
  aspect ratios $\epsilon$ as indicated in the plot.  The number of
  iterations required to converge increases for decreasing $\epsilon$, 
  with a long ``plateau'' for small aspect ratios.}
\end{figure}

\begin{figure}
\begin{center}
\includegraphics*[width=12cm]{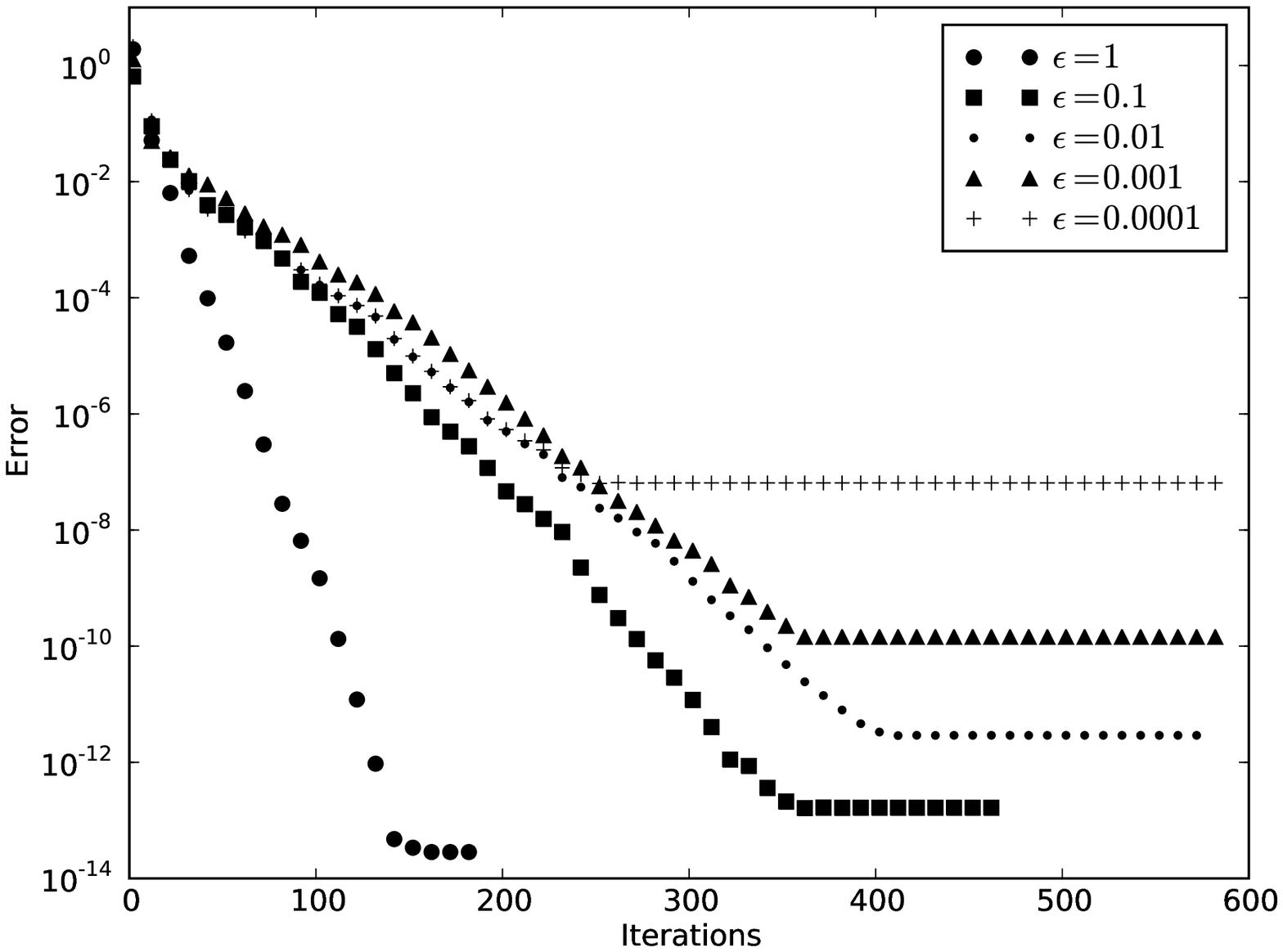}
\end{center}
\caption{\label{monster_vl_eps} Plot showing error (using the
  $\infty$-norm) against number of iterations, for the Conjugate
  Gradient method applied to the Poisson equation discretised on mesh
  B (shown in figure \ref{monster}, $\epsilon=0.001$), using the vertically lumped
  preconditioner. The mesh has been rescaled to various different
  aspect ratios $\epsilon$ as indicated in the plot.  The
  convergence rate becomes independent of $\epsilon$ for small
  aspect ratios.}
\end{figure}

\begin{figure}
\begin{center}
\includegraphics*[width=12cm]{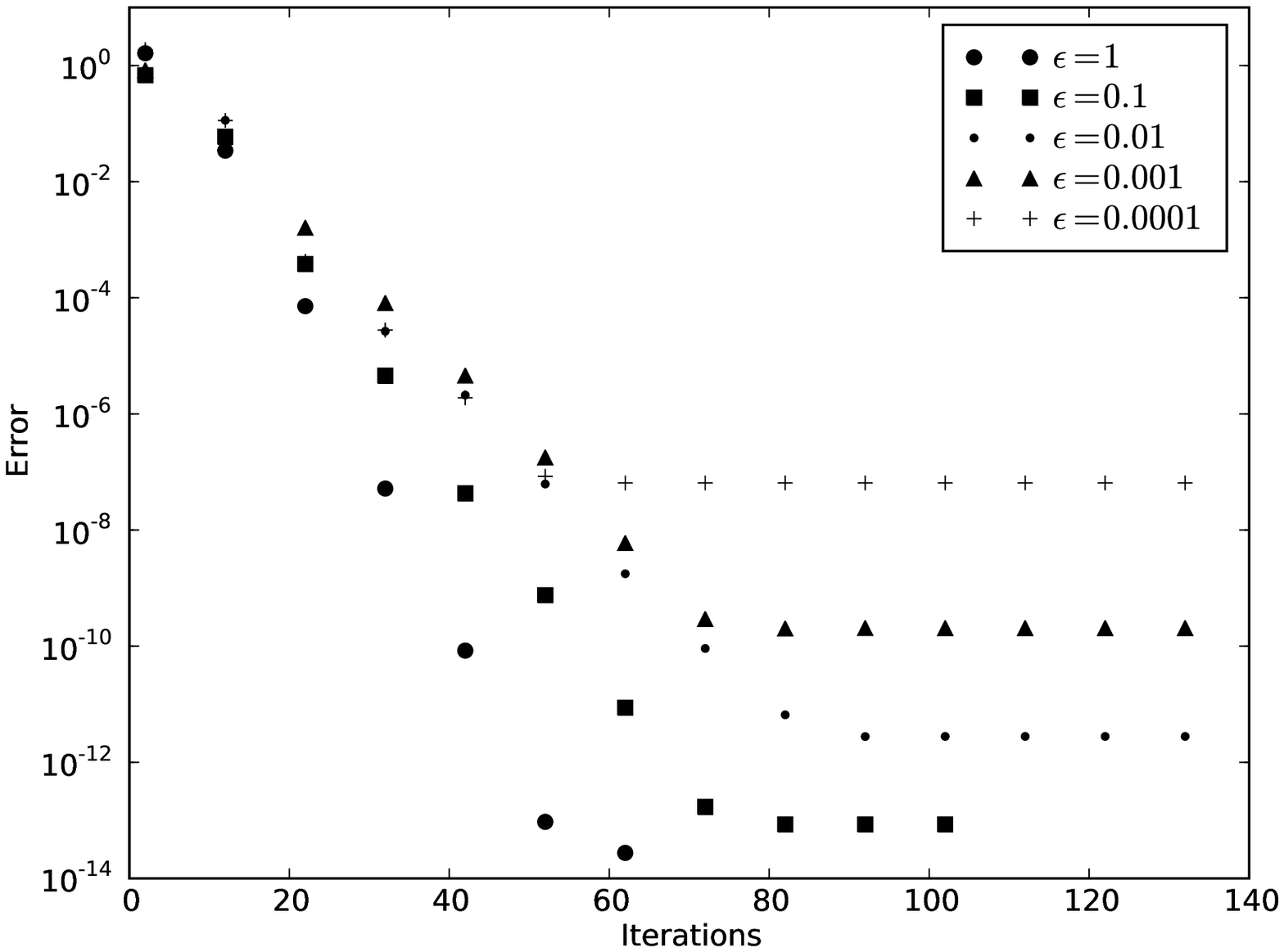}
\end{center}
\caption{\label{monster_vl_as_eps} Plot showing error (using the
  $\infty$-norm) against number of iterations, for the Conjugate
  Gradient method applied to the Poisson equation discretised on mesh
  B (shown in figure \ref{monster}), using the vertically lumped
  preconditioner combined with additive smoothing. The mesh has been
  rescaled to various different aspect ratios $\epsilon$ as
  indicated in the plot. The convergence rate becomes independent of
  $\epsilon$ for small aspect ratios, with a substantial improvement
  over the vertically lumped preconditioner with standard SOR smoothing,
  shown in figure \ref{monster_vl_eps}.}
\end{figure}

Next we present results for Mesh B which is a multiscale mesh, as
illustrated in figure \ref{monster}. In figure
\ref{monster_no_vl_eps}, we plot the convergence of the smoothed
aggregation multigrid method applied to the matrix obtained from mesh
B, which again shows a convergence plateau which becomes longer as
$\epsilon\to 0$. Figures \ref{monster_vl_eps} and
\ref{monster_vl_as_eps} show the convergence rate in iterations for
the vertically lumped preconditioner with and without the additive
smoother respectively. In both cases the convergence rate becomes
independent of $\epsilon$ for small aspect ratios. In this case the
additive smoother is producing a faster decay of error as the number
of iterations increases, compared to the SOR smoother. This suggests
that the additive smoother is more effective at approximating the
baroclinic components of the solution, which have a complex multiscale
structure. The additive smoother uses an algebraic multigrid cycle
applied to the baroclinic components, which operates at several scales
simultaneously.

\begin{figure}
\begin{center}
\includegraphics*[width=12cm]{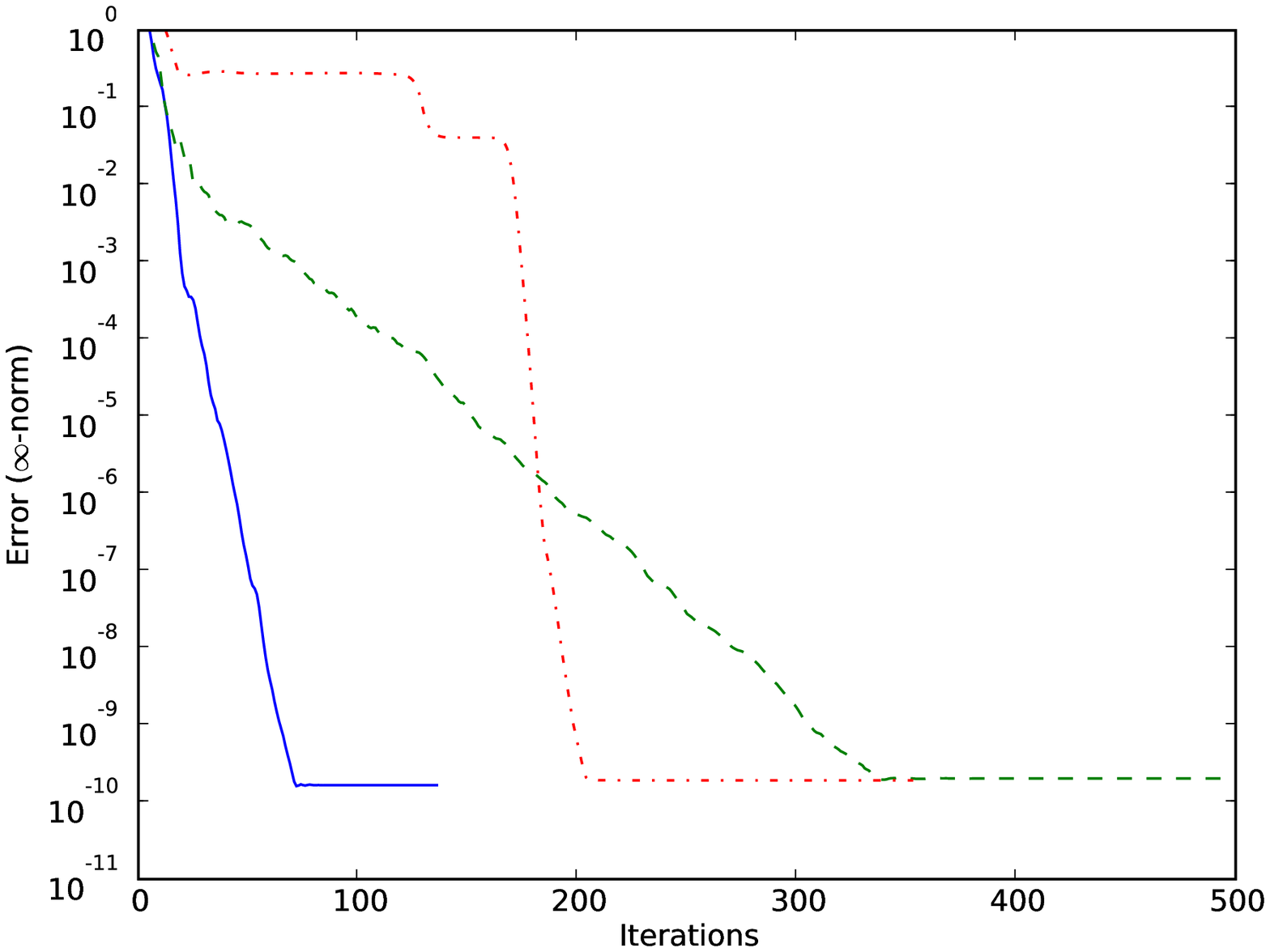}
\end{center}
\caption{\label{eps0_001_its} Plot showing error (using the
  $\infty$-norm) against number of iterations, for the Conjugate
  Gradient method applied to the Poisson equation discretised on mesh
  B (shown in figure \ref{monster}), with various different
  preconditioners. The continuous line indicates the vertically lumped
  preconditioner with the additive smoother, the dashed line indicates
  the vertically lumped preconditioner without the additive smoother,
  and the dash dotted line indicates the smoothed aggregation multigrid
  preconditioner.}
\end{figure}

Finally in figure \ref{eps0_001_its}, the error is plotted against number of
iterations for the smoothed aggregation preconditioner and the
vertically lumped preconditioner with and without the additive
smoother for the matrix obtained from Mesh B with $\epsilon=0.001$.
We observe again that the vertically lumped preconditioner converges
much faster than the smoothed aggregation preconditioner. Here the
vertically lumped preconditioner does not exhibit a convergence
plateau but does have a slow rate of convergence which we attribute to
the presence of the small eigenvalues associated with small scales in
the mesh which are not altered by the vertically lumped
preconditioner.  The inclusion of the additive smoother means that the
number of iterations is dramatically reduced, since the additive
smoother is a multigrid preconditioner which treats all of the scales
in the mesh. Despite the added cost of the additive smoother, we
observe that it results in a much more efficient solver. We conclude
that the additive smoother should be used when small scales are
present in the mesh which lead to eigenvalues which are of the same
size as those associated with eigenvectors corresponding to horizontal
modes.

\begin{figure}
\begin{center}
\includegraphics*[width=12cm]{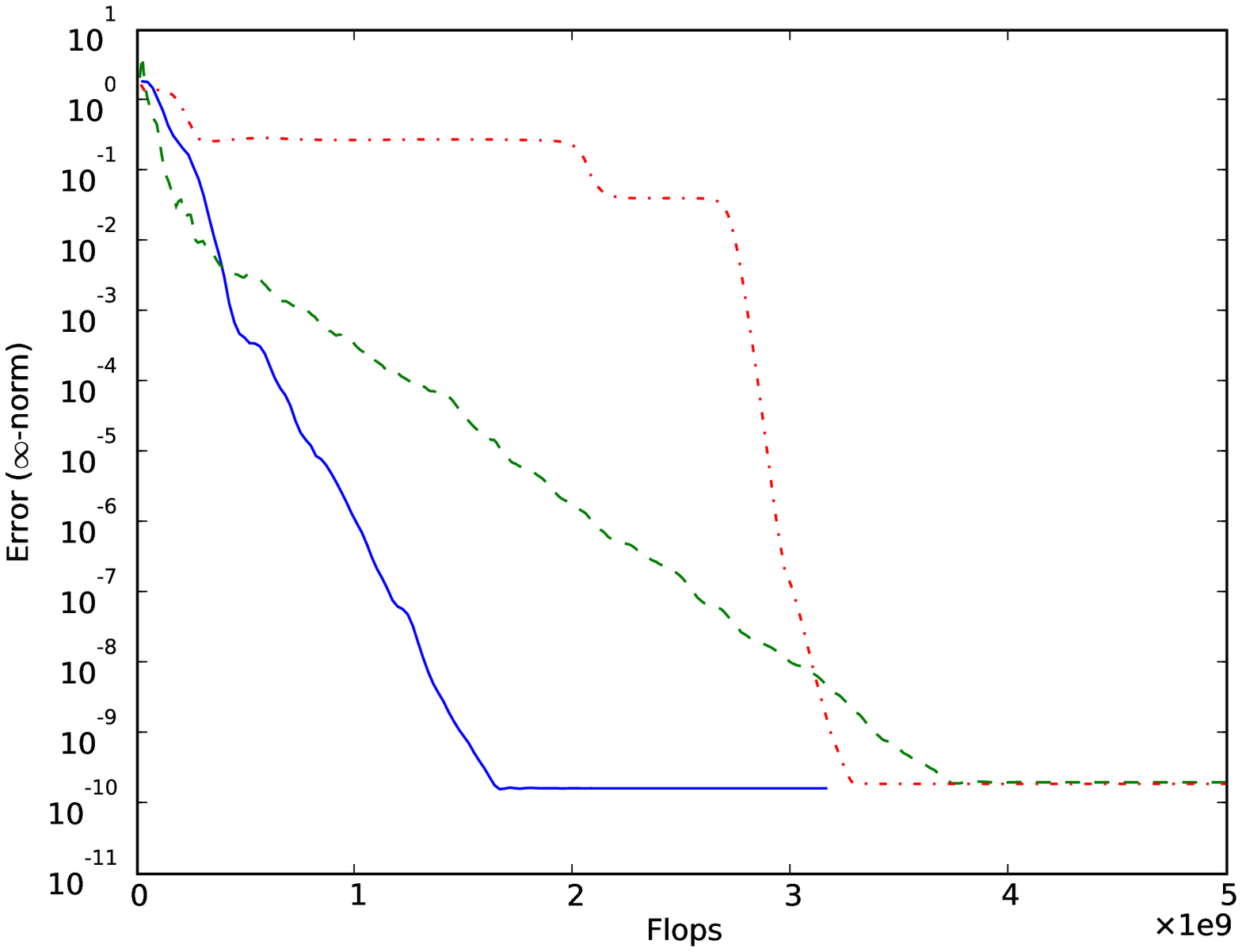}
\end{center}
\caption{\label{eps0_001_flops} Plot showing error (using the
  $\infty$-norm) against number of floating point operations (flops)
  counted by PETSc, for the Conjugate Gradient method applied
  to the Poisson equation discretised on Mesh B (shown in figure
  \ref{monster}), with various different preconditioners. The
  continuous line indicates the vertically lumped preconditioner with
  the additive smoother, the dashed line indicates the
  vertically lumped preconditioner without the additive smoother, and
  the dash dotted line indicates the smoothed aggregation multigrid
  preconditioner.}
\end{figure}

\section{Summary and outlook}
\label{summary}

In this paper we discussed the ill conditioning of the linear system
obtained from the finite element approximation of the pressure Poisson
equation on general vertically unstructured meshes in small aspect
ratio domains (such as the global oceans). We showed that the
condition number scales like $\epsilon^{-2}$ as $\epsilon \to 0$
(where $\epsilon$ is the aspect ratio $H/L$) in the case in which Neumann
boundary conditions are set on all surfaces. We also showed that the
condition number is independent of $\epsilon$ when Dirichlet
conditions are applied at the top surface.  This motivated a
preconditioner consisting of two stages: in the first stage an
approximate reduced system for the surface degrees of freedom is
solved, and in the second stage the solution is reconstructed
throughout the domain with the approximate surface solution used as a
Dirichlet boundary condition. The first stage results in a much
smaller linear system, and the second stage involves a submatrix which
has a condition number which is independent of $\epsilon$. The reduced
system is obtained using an algebraic multigrid prolongation operator
which approximates the vertical extrapolation operator, and the second
stage submatrix can be solved using further algebraic multigrid
stages.  Using numerical experiments, we showed that this
preconditioner, when combined with the Conjugate Gradient method,
results in a solver which has a convergence rate that is independent
of the aspect ratio. Further, we showed that the additional
computational cost of using the additive smoother means that it is
only beneficial in truly multiscale meshes. Those are meshes 
that do not just have two entirely different lenght scales, 
the horizontal and the vertical, but a whole range of scales inbetween
. This strategy will
become crucial when solving process study problems consisting of small
scale dynamics (such as open ocean deep convection, or density
overflows) that are embedded in a large scale circulation.  We also
anticipate that the smoother will become important when parallel
domain decomposition methods are used, where (block) SOR smoothing
methods are known to be less effective. The methods described in this
paper have been implemented in ICOM where they will be used to
investigate adaptive unstructured large scale ocean modelling.

\section{Acknowledgements}
The authors would like to acknowledge the funding of the UK 
Natural Environment Research Council under grant NE/C52101X/1, 
and of Fujitsu Laboratories of Europe.

\bibliographystyle{elsarticle-harv}
\bibliography{pressure}

\end{document}